\documentclass[a4paper,12pt]{article}
\begin{document}

\title{Integral operator approach over octonions to solution of nonlinear PDE.}
\author{E. Frenod, S.V. Ludkowski}

\date{05 May 2016}
\maketitle

\begin{abstract}
Integration of nonlinear partial differential equations with the
help of the non-commutative integration over octonions is studied.
An apparatus permitting to take into account symmetry properties of
PDOs is developed. For this purpose formulas for calculations of
commutators of integral and partial differential operators are
deduced. Transformations of partial differential operators and
solutions of partial differential equations are investigated.
Theorems providing solutions of nonlinear PDEs are proved. Examples
are given. Applications to PDEs of hydrodynamics and other types
PDEs are described.

\footnote{key words and phrases: hypercomplex, octonion algebra, nonlinear
partial differential equation, non-commutative integration,
integral operator 
\\ addresses:
E. Frenod, LMBA, UMR CNRS 6205,
Laboratoire de Math\'ematiques de Bretagne Atlantique,
Universit\'e de Bretagne Sud, Campus de Tohannic BP573
-   56017 Vannes, France \\
Emmanuel.Frenod@univ-ubs.fr \\
S.V. Ludkovsky, LMBA, UMR CNRS 6205, Laboratoire de Math\'ematiques
de Bretagne Atlantique, Universit\'e de Bretagne Sud, Campus de
Tohannic BP573 -   56017 Vannes, France \par and Department of
Applied Mathematics, Moscow State Technical University MIREA, av.
Vernadsky 78, Moscow
119454, Russia\\
Ludkowski@mirea.ru }

\end{abstract}

\section{Introduction.}
Analysis over hypercomplex numbers develops fast and has important
applications in geometry and partial differential equations
including that of nonlinear (see \cite{br} - \cite{guesprqa},
\cite{lucveleq2013}-\cite{lufjms2010} and references therein). As a
consequence it gives new opportunities for integration of different
types of partial differential equations (PDEs).
 It is worth to mention that the quaternion skew field ${\bf H}={\cal
 A}_2$, the octonion algebra ${\bf O}={\cal A}_3$ and Cayley-Dickson
algebras ${\cal A}_r$ have found a lot of applications not only in
mathematics, but also in theoretical physics (see \cite{br} -
\cite{guespr} and references therein).
\par Mixed type PDEs play very important role not only in mathematics,
but also in physics. For example, they describe two-dimensional
motions of a stratified rotating liquid, electromagnetic fields in
crystals, internal gravitational waves, non stationary filtration
process of liquid in a fissure porous medium, dissipation process,
cold plasma, two temperature plasma in an external magnetic field,
etc.  (see \cite{svalkorplb} and references therein). They were
studied with the real time and two spatial variables. But there are
needs to integrate more general PDE of such type with larger number
of variables or their systems.
\par Frequently types of considered PDE in modeling different processes and
in studying their solutions are restricted by available tools of
mathematical analysis. Otherwise numerical methods and computer
calculations are used. But in many cases it is also necessary to
analyze properties of solutions. Therefore analytic approaches apart
from that of numerical provide in this respect many advantages. On
the other hand, ranges of mathematical analysis strongly depend on
used number systems. Classically in mathematical analysis and PDEs
real and complex fields are used. But hypercomplex analysis enlarges
its scopes. \par This article is devoted to analytic approaches to
solution of PDEs and taking into account their symmetry properties.
For this purpose the octonion algebra is used. This is actual
especially in recent period because  of increasing interest to
non-commutative analysis and its applications. It is worth to
mention that each problem of PDE can be reformulated using the
octonion algebra. The approach over octonions enlarges a class of
PDEs which can be analytically integrated in comparison with
approaches over the real field and the complex field.
\par  There is stimulus for investigations caused by needs to integrate known PDEs
(see, for example, \cite{hormb3v,hormbl,kneschkeb,svalkorplb} and
references therein) and by the progress of algebra
\cite{baez,kurosh,serodio,schaferb}. Certainly algebras are widely
used in PDEs (see also, for example, \cite{emch,guetze,
prastaro1}-\cite{prastaro3} and references therein). But in previous
works mainly associative algebras were used for such purposes.
\par On the other hand nonlinear PDEs frequently are more complicated
and demand specific approaches to get their solutions
\cite{hormbl,kichenb96,polzayrb,svalkorplb} in comparison with that
of linear.
\par We exploit a new approach based on the non-commutative
integration over non-associative Cayley-Dickson algebras that to
integrate definite types of nonlinear PDEs. This work develops
further results of the previous article \cite{lucveleq2013}. The
obtained below results open new perspectives and permit to integrate
nonlinear PDEs with variable coefficients and analyze symmetries of
solutions as well.
\par In the following sections integration of nonlinear PDEs with the
help of the non-commutative integration over quaternions, octonions
and Cayley-Dickson algebras is studied. For this purpose formulas
for calculations of commutators of integral and partial differential
operators are deduced. Transformations of partial differential
operators and solutions of partial differential equations are
investigated. An apparatus permitting to take into account symmetry
properties of PDOs is developed. Theorems providing solutions of
nonlinear PDEs are proved. Examples are given. Applications to PDEs
used in hydrodynamics and other types PDEs are described. The
results of this paper can be applied to integration of some kinds of
nonlinear Sobolev type PDEs as well.
\par  All main results of this paper are obtained
for the first time. They can be used for further investigations of
PDEs and properties of their solutions. For example, generalized
PDEs including terms such as ${\Delta }^p$ or ${\nabla }^p$ for
$p>0$ or even complex $p$ can be investigated.

\section{Integral operators over octonions.}
\par To avoid misunderstandings we first present our definitions and notations.
\par {\bf 1. Notations and Definitions.}
By ${\cal A}_r$ we denote the Cayley-Dickson algebra over the real field $\bf R$ with
generators $i_0,...,i_{2^r-1}$ so that $i_0=1$, $~i_j^2=-1$ for
each $j\ge 1$, $~i_ji_k=-i_ki_j$ for each $j\ne k\ge 1$, $2\le r\in \bf N$.
\par Henceforward PDEs are considered on a domain $U$ in ${\cal A}_r^m$ such that
\par $(D1)$ each projection ${\bf p}_j(U)=:U_j$ is
$(2^r-1)$-connected;
\par $(D2)$ $\pi _{{\sf s},{\sf p},{\sf t}}(U_j)$ is simply
connected in $\bf C$ for each $k=0,1,...,2^{r-1}$, ${\sf s}=i_{2k}$,
${\sf p}=i_{2k+1}$, ${\sf t}\in {\cal A}_{r,{\sf s},{\sf p}}$ and
${\sf u}\in {\bf C}_{{\sf s},{\sf p}}$, for which
a Cayley-Dickson number $z$ exists satisfying the condition $z={\sf u}+{\sf t}\in U_j$, \\
where $e_j = (0,...,0,1,0,...,0)\in {\cal A}_r^m$ is the vector with
$1$ on the $j$-th place, ${\bf p}_j(z) = \mbox{ }^jz$ for each $z\in
{\cal A}_r^m$, where the decomposition is used $z=\sum_{j=1}^m\mbox{ }^jz e_j$, $\mbox{ }^jz\in
{\cal A}_r$ for each $j=1,...,m$, $m\in {\bf N} := \{ 1,2,3,... \}
$,  the projections are the following $~\pi _{{\sf s},{\sf p},{\sf t}}(V):= \{ {\sf u}: z\in V,
z=\sum_{{\sf v}\in \bf b}w_{\sf v}{\sf v},$ ${\sf u}=w_{\sf s}{\sf
s}+w_{\sf p}{\sf p} \} $ for a domain $V$ in ${\cal A}_r$ for each
${\sf s}\ne {\sf p}\in \bf b$, where ${\sf t}:=\sum_{{\sf v}\in {\bf
b}\setminus \{ {\sf s}, {\sf p} \} } w_{\sf v}{\sf v} \in {\cal
A}_{r,{\sf s},{\sf p}}:= \{ z\in {\cal A}_r:$ $z=\sum_{{\sf v}\in
\bf b} w_{\sf v}{\sf v},$ $w_{\sf s}=w_{\sf p}=0 ,$ $w_{\sf v}\in
\bf R$ $\forall {\sf v}\in {\bf b} \} $, where ${\bf b} := \{
i_0,i_1,...,i_{2^r-1} \} $ is the family of standard generators of
the Cayley-Dickson algebra ${\cal A}_r$. Frequently we take $m=1$.
Henceforth, we consider a domain $U$ satisfying Conditions
$(D1,D2)$ if something other is not outlined.

\par {\bf 2. Operators.} An ${\bf R}$ linear space $X$ which is also
left and right ${\cal A}_r$ module will be called an ${\cal A}_r$
vector space. It is supposed that $X$ can be presented as the direct
sum $X=X_0i_0\oplus ... \oplus X_{2^r-1} i_{2^r-1}$, where $X_0$,...,$X_{2^r-1}$ are
pairwise isomorphic real linear spaces. Particularly, for $r=2$ this
module is associative: $(xa)b=x(ab)$ and  $(ab)x=a(bx)$ for all
$x\in X$ and $a,b \in {\bf H}$, since the quaternion skew field
${\cal A}_2={\bf H}$ is associative. This module is alternative for
$r=3$: $(xa)a=x(a^2)$ and $(a^2)x=a(ax)$ for all $x\in X$ and $a\in
{\bf O}$, since the octonion algebra ${\bf O}={\cal A}_3$ is
alternative. \par  Let $X$ and $Y$ be two $\bf R$ linear normed
spaces which are also left and right ${\cal A}_r$ modules, where
$2\le r$, such that \par $(1)$ $0\le \| ax \|_X \le |a| \|x \|_X $
and $ \| x a \|_X \le |a| \|x \|_X$ for all $x\in X$ and $a\in {\cal A}_r$
and \par $(2)$ $\| x+y \|_X \le \| x \|_X + \| y \|_X $ for all $x, y\in X$ and
and \par $(3)$ $\| bx \|_X = |b| \|x \|_X = \| x b \|_X$ for each $b\in \bf R$
and $x\in X$, where for $r=2$ and $r=3$ Condition $(1)$ takes the form
\par $(1')$ $0\le \| ax \|_X = |a| \|x \|_X = \| x a \|_X$ for all $x\in X$ and $a\in {\cal A}_r$.
\par Such spaces $X$ and $Y$ will be called ${\cal A}_r$ normed spaces.
\par An ${\cal A}_r$ normed space complete
relative to its norm will be called an ${\cal A}_r$ Banach space.
\par Put $X^{\otimes k} := X\otimes _{\bf R} ... \otimes
_{\bf R} X$ to be the $k$ times ordered tensor product over $\bf R$
of $X$.  By $L_{q,k}(X^{\otimes k},Y)$ we denote a family of all
continuous $k$ times $\bf R$ poly-linear and ${\cal A}_r$ additive
operators from $X^{\otimes k}$ into $Y$. If $X$ and $Y$ are normed
${\cal A}_r$ spaces and $Y$ is complete relative to its norm, then
$L_{q,k}(X^{\otimes k},Y)$ is also a normed $\bf R$ linear and left
and right ${\cal A}_r$ module complete relative to its norm. In
particular, $L_{q,1}(X,Y)$ is denoted by $L_q(X,Y)$ as well.
\par  If $A\in L_q(X,Y)$ and $A(xb)=(Ax)b$
or $A(bx)=b(Ax)$ for each $x\in X_0$ and $b\in {\cal A}_r$, then an
operator $A$ we call right or left ${\cal A}_r$-linear respectively.
\par An $\bf R$ linear space of left (or right) $k$ times ${\cal A}_r$
poly-linear operators is denoted by $L_{l,k}(X^{\otimes k},Y)$ (or
$L_{r,k}(X^{\otimes k},Y)$ respectively).
\par An $\bf R$-linear operator $A: X\to X$ will be called right or
left strongly ${\cal A}_r$ linear if
\par $A(xb)=(Ax)b$ or
\par $A(bx)=b (Ax)$ for each $x\in X$ and $b\in {\cal A}_r$ correspondingly.
\par An ${\bf R}$ linear ${\cal A}_r$ additive operator $A$ is called
invertible if it is densely defined and one-to-one and has a dense
range ${\cal R}(A)$.
\par Henceforward, if an expression of the form
\par $(4)$ $\sum_k [(I-A_x)~\mbox{}_kf(x,y)]~\mbox{}_kg(y)=u(x,y)$ \\ will appear on a domain $U$,
which need to be inverted we consider the case when \par $(RS)$
$(I-A_x)$ is either right strongly ${\cal A}_r$ linear, or right
${\cal A}_r$ linear and $~\mbox{}_kf\in X_0$ for each $k$, or ${\bf
R}$ linear and $~\mbox{}_kg(y)\in {\bf R}$ for each $k$ and every
$y\in U$, at each point $x\in U$, since $\bf R$ is the center
of the Cayley-Dickson algebra ${\cal A}_r$, where $2\le r$.

\par {\bf 3. First order PDOs.} We consider an arbitrary first order partial
differential operator $\sigma $ given by the formula
\par $(1)$ $\sigma f= \sum_{j=0}^{2^r-1}
i_j^* (\partial f/\partial z_{\xi (j)}) {\psi }_j$,
\\ where $f$ is a differentiable ${\cal A}_r$-valued function on the domain $U$ satisfying
Conditions 1$(D1,D2)$, $~ 2\le r$, $ ~ i_0,...,i_{2^r-1}$ are the
standard generators of the Cayley-Dickson algebra ${\cal A}_r$, $a^*={\tilde a} := a_0i_0-a_1i_1-...
-a_{2^r-1} i_{2^r-1}$ for each $a = a_0i_0+a_1i_1+... +a_{2^r-1} i_{2^r-1}$ in ${\cal A}_r$
with $a_0,...,a_{2^r-1}\in \bf R$; $~\psi _j$ are real constants so that $\sum_j\psi _j^2> 0$, $ ~ \xi :
\{ 0,1,...,2^r-1 \} \to \{ 0,1,...,2^r-1 \}$ is a surjective
bijective mapping, i.e. $\xi $ belongs to the symmetric group
$S_{2^r}$ (see also \S 2 in \cite{lucmft12}).
\par For an ordered product $\{ \mbox{}_1f...\mbox{}_kf \}_{q(k)}$
of differentiable functions $\mbox{}_sf$ we put
\par $(2)$ $\mbox{}^s\sigma \{ \mbox{}_1f...\mbox{}_kf \}_{q(k)}
= \sum_{j=0}^n  i_j^* \{ \mbox{}_1f...(\partial \mbox{}_sf/\partial
z_{\xi (j)})...\mbox{}_kf \} _{q(k)} \psi _j,$ \\
where a vector $q(k)$ indicates on an order of the multiplication in
the curled brackets (see also \S  2 \cite{ludoyst,ludfov}), so that
\par $(3)$ $\sigma \{ \mbox{}_1f...\mbox{}_kf \}_{q(k)}=\sum_{s=1}^k
\mbox{}^s\sigma \{ \mbox{}_1f...\mbox{}_kf \}_{q(k)}$. \par
\par Symmetrically other operators
\par $(4)$ ${\hat \sigma }f= \sum_{j=0}^{2^r-1}
(\partial f/\partial z_{\xi (j)}) i_j {\psi }_j$, \\
are defined. Therefore, these operators are related by the formula:
\par $(5)$ $(\sigma f)^* = {\hat \sigma } (f^*) $.
\par Operators $\sigma $ given by $(1)$ are right ${\cal A}_r$
linear.

\par {\bf 4. Integral operators.} We consider integral operators of the form:
$$(1)\quad K(x,y) = F(x,y) +p\mbox{}_{\sigma } \int_x^{\infty }
F(z,y)N(x,z,y)dz ,$$ where $\sigma $ is an $\bf
R$-linear partial differential operator as in \S 3 and
$\mbox{}_{\sigma } \int $ is the non-commutative line integral
(anti-derivative operator) over the Cayley-Dickson algebra ${\cal
A}_r$ from \cite{lucmft12} or \S 4.2.5
\cite{ludancdnb}, where $F$ and $K$ are continuous
functions with values in the Cayley-Dickson algebra ${\cal A}_r$ or
more generally in the real algebra $Mat_{n\times
n}({\cal A}_r)$ of $n\times n$ matrices with entries in ${\cal A}_r$,
$p$ is a nonzero real parameter.
For definiteness we take the right ${\cal A}_r$ linear anti-derivative operator
$\mbox{}_{\sigma } \int g(z)dz$.
\par Let a domain $U$ be provided with a foliation
by locally rectifiable paths $\{ \gamma ^{\alpha }: ~ \alpha \in
\Lambda \} $ (see also \cite{lucmft12} or
\cite{ludancdnb}), where $\Lambda $ is a set described below. We take for
definiteness a canonical closed domain $U$ in ${\hat {\cal A}}_r$
satisfying Conditions 1$(D1,D2)$ so that $\infty \in U$, where
${\hat {\cal A}}_r = {\cal A}_r \cup \{ \infty \} $ denotes the
one-point compactification of ${\cal A}_r$, $ ~ 2\le r<\infty $.
\par A domain $U$ is called foliated by locally rectifiable paths $\{ \gamma
^{\alpha }: ~ \alpha \in \Lambda \} $ if $\gamma : <a_{\alpha
},b_{\alpha }> \to U$ for each $\alpha $ and it satisfies the
following three conditions: \par $(F1)$ $\bigcup_{\alpha \in \Lambda
} \gamma ^{\alpha }(<a_{\alpha },b_{\alpha }>) =U$ and \par $(F2)$
$\gamma ^{\alpha }(<a_{\alpha },b_{\alpha }>)\cap \gamma ^{\beta
}(<a_{\beta },b_{\beta }>) = \emptyset $ for each $\alpha \ne \beta
\in \Lambda $. \\ Moreover, if the boundary $\partial U = cl
(U)\setminus Int (U)$ of the domain $U$ is non-void then \par $(F3)$
$\partial U = (\bigcup_{\alpha \in \Lambda _1} \gamma ^{\alpha
}(a_{\alpha }))\cup
(\bigcup_{\beta \in \Lambda _2} \gamma ^{\beta }(b_{\beta }))$, \\
where $\Lambda _1 = \{ \alpha \in \Lambda : <a_{\alpha },b_{\beta
}>= [a_{\alpha },b_{\beta }> \} $, $\Lambda _2 = \{ \alpha \in
\Lambda : <a_{\alpha },b_{\beta }>= <a_{\alpha },b_{\beta }] \} $.
For the canonical closed subset $U$ we have $cl (U)=U=cl (Int(U))$,
where $cl (U) $ denotes the closure of $U$ in ${\cal A}_v$ and $Int
(U)$ denotes the interior of $U$ in ${\cal A}_v$. For convenience
one can choose a $C^n$ foliation, that is each $\gamma ^{\alpha }$ is of
class $C^n$, where $n\in \bf N$. When $U$ is with non-void boundary we choose a
foliation family such that $\bigcup_{\alpha \in \Lambda } \gamma
(a_{\alpha }) =\partial U_1$, where a set $\partial U_1$ is open in
the boundary $\partial U$ and so that $w|_{\partial U_1}$ would be a
sufficient initial condition to characterize a unique branch of an
anti-derivative $w(x)={\cal I}_{\sigma }f(x) = ~ \mbox{}_{\sigma }
\int_{\mbox{}_0x}^x f(z)dz$, \\ where $\mbox{}_0x\in \partial U_1$,
$x\in U$, $\gamma ^{\alpha }(t_0)=\mbox{}_0x$, $\gamma ^{\alpha }
(t)=x$ for some $\alpha \in \Lambda $, $t_0$ and $t\in <a_{\alpha
},b_{\alpha }>$, \par $~ \mbox{}_{\sigma } \int_{\mbox{}_0x}^x
f(z)dz= ~ \mbox{}_{\sigma } \int_{\gamma ^{\alpha }|_{[t_0,t]}}
f(z)dz$.
\par In accordance with Theorems 2.4.1 and 2.5.2 \cite{lucmft12} or
4.2.5 and 4.2.23 \cite{ludancdnb} the equality
\par $(2)$ $\sigma _x ~ \mbox{}_{\sigma } \int_{\mbox{}_0x}^x
g(z)dz =g(x)$ \\
is satisfied for a continuous function $g$ on a domain $U$ as in \S
1 and a foliation as above.
\par Particularly in the class of ${\cal A}_r$ holomorphic functions
in the domain satisfying Conditions 1$(D1,D2)$ this line integral
depends only on initial and final points due to the homotopy theorem
\cite{ludoyst,ludfov}.
\par  We denote by ${\bf P}={\bf P}(U)$ the family of all locally rectifiable paths $\gamma :
<a_{\gamma },b_{\gamma }>\to U$ supplied with the family of
pseudo-metrics
\par $(3)$ $\rho ^{a,b,c,d}(\gamma ,\omega ) := |\gamma (a)-\omega (c)|+
\inf_{\phi } V_a^b(\gamma (t) - \omega (\phi (t))$ \\
where the infimum is taken by all diffeomorphisms $\phi : [a,b]\to
[c,d]$ so that $\phi (a)=c$ and $\phi (b)=d$, $~ a<b$, $~ c<d$,
$[a,b]\subset <a_{\gamma }, b_{\gamma }>$, $[c,d]\subset <a_{\omega
}, b_{\omega }>$. We take a foliation such that $\Lambda $ is a
uniform space and  the limit
\par $(4)$ $\lim_{\beta \to \alpha }\rho ^{a,b,a,b}(\gamma ^{\beta },\gamma ^{\alpha
})=0$ \\ is zero for each $[a,b]\subset <a_{\alpha }, b_{\alpha }>$.
\par  For example, we can take $\Lambda = {\bf R}^{n-1}$ for a foliation of the entire
Cayley-Dickson algebra ${\cal A}_r$, where $n=2^r$, $t\in {\bf R}$,
$\alpha \in \Lambda $, so that $\bigcup_{\alpha \in \Lambda } \gamma
^{\alpha }([0,\infty )) = \{ z\in {\cal A}_r: ~ Re((z-y)v^*)\ge 0 \}
$ and $\bigcup_{\alpha \in \Lambda } \gamma ^{\alpha }((-\infty ,0])
= \{ z\in {\cal A}_r: ~ Re((z-y)v^*)\le 0 \} $ are two real
half-spaces, where $v, y\in {\cal A}_r$ are marked Cayley-Dickson
numbers and $v\ne 0$. Particularly, we can choose the foliation such
that $\gamma ^{\alpha } (0)=y+\alpha _1v_1+...+\alpha
_{2^r-1}v_{2^r-1}$ and $\gamma ^{\alpha }(t)=tv_0+\gamma ^{\alpha
}(0)$ for each $t\in {\bf R}$, where $v_0,...,v_{2^r-1}$ are $\bf
R$-linearly independent vectors in ${\cal A}_r$.

\par Therefore, the expression $$(5)\quad \mbox{}_{\sigma } \int_x^{\infty }
g(z)dz := \mbox{}_{\sigma } \int_{\gamma ^{\alpha }|_{[t_x,b_{\alpha
}>}} g(z)dz $$ denotes a non-commutative line integral over ${\cal
A}_r$ along a path $\gamma ^{\alpha }$ so that $\gamma ^{\alpha
}(t_x)=x$ and $\lim_{b\to b_{\alpha }}\gamma ^{\alpha }(t)=\infty $
for an integrable function $g$, where $t_x\in <a_{\alpha },b_{\alpha
}>$, $ ~ \alpha \in \Lambda $, $a_{\alpha }=a_{\gamma ^{\alpha }}$,
$b_{\alpha }=b_{\gamma ^{\alpha }}$. It is sometimes convenient to
use the line integral
$$(6)\quad \mbox{}_{\sigma } \int_x^{-\infty }
g(z)dz := \mbox{}_{\sigma } \int_{\gamma ^{\alpha }|_{<a_{\alpha
},t_x]}} g(z)dz ,$$ when $\lim_{a\to a_{\alpha }}\gamma ^{\alpha
}(t)= \infty $.
\par Put for convenience $\sigma ^0=I$, where $I$ denotes the unit
operator, $ ~ \sigma ^m$ denotes the $m$-th power of $\sigma $ for
each non-negative integer $0\le m\in {\bf Z}$.

\par {\bf 5. Proposition.} {\it Let $F\in C^m(U^2,Mat_{n\times n}({\cal A}_r))$
and $N\in C^m(U^3,Mat_{n\times n}({\cal A}_r))$ and let
$$(1)\quad \lim_{z\to \infty } ~ \mbox{}^1\sigma _z^k ~ \mbox{}^2\sigma
_x^s ~ \mbox{}^2\sigma _z^lF(z,y)N(x,z,y)=0$$ for each $x,
y$ in a domain $U$ satisfying Conditions 1$(D1,D2)$ with $\infty
\in U$ and every non-negative integers $0\le k, s, l\in {\bf Z}$
such that $k+s+l\le m $. Suppose also that $\mbox{}_{\sigma }
\int_x^{\infty } \partial ^{\alpha }_x\partial ^{\beta }_y \partial
^{\omega }_z [F(z,y) N(x,z,y)] dz $ converges uniformly by
parameters $x, y $ on each compact subset $W\subset U\subset {\cal
A}_r^2$ for each $|\alpha |+|\beta |+|\omega |\le m$, where $\alpha
=(\alpha _0,...,\alpha _{2^r-1})$, $|\alpha |=\alpha _0+...+\alpha
_{2^r-1}$, $\partial ^{\alpha }_x=\partial ^{|\alpha |}/\partial
x_0^{\alpha _0} ...\partial x_{2^r-1}^{\alpha _{2^r-1}}$. Then the
non-commutative line integral $\mbox{}_{\sigma } \int_x^{\infty }
F(z,y)N(x,z,y)dz $ from \S 4 satisfies the identities:
$$(2)\quad \sigma ^m_x ~ \mbox{}_{\sigma}\int_x^{\infty } F(z,y)N(x,z,y)dz =
~ \mbox{}^2\sigma ^m_x ~ \mbox{}_{\sigma}\int_x^{\infty } F(z,y) N(x,z,y)dz + A_m(F,N)(x,y),$$
$$(3)\quad \mbox{}^1\sigma ^m_z~ \mbox{}_{\sigma}\int_x^{\infty }  F(z,y) N(x,z,y)dz =
(-1)^m ~ \mbox{}^2\sigma ^m_z ~\mbox{}_{\sigma}\int_x^{\infty }
F(z,y)N(x,z,y)dz + B_m(F,N)(x,y),$$ where
$$(4)\quad A_m(F,N)(x,y)= - ~\mbox{}^2\sigma
^{m-1}_x [F(x,y)N(x,z,y)]|_{z=x} + \sigma _x~ A_{m-1}(
F,N)(x,y)$$ for $m\ge 2$, $$(5) \quad B_m(F(z,y),
N(x,z,y)) = (-1)^m ~ \mbox{}^2\sigma ^{m-1}_z F(x,y)
N(x,z,y) + \mbox{}^1\sigma _z~ B_{m-1}(F(z,y),
N(x,z,y))$$ for $m\ge 2$, $B_m(F,N)(x,y)=B_m(F(z,y),
N(x,z,y))|_{z=x}$;
\par $(6)\quad  ~ A_1(F,N)(x,y)= - F(x,y)N(x,x,y)$,
\par $(7)$ $B_1(F(z,y),N(x,z,y))= - F(z,y)N(x,z,y)$, \\
 $ ~ \sigma _x$ is an operator $\sigma $ acting by the
variable $x\in U\subset {\cal A}_r$.}
\par {\bf Proof.} Using the conditions of this proposition and the
theorem about differentiability of improper integrals by parameters
(see, for example, Part IV, Chapter 2, \S 4 in \cite{kamyn})
we get the equality \par $\mbox{}_{\sigma } \int_x^{\infty
}\partial ^{\alpha }_x\partial ^{\beta }_y  \partial ^{\omega }_z
[F(z,y) N(x,z,y)] dz = \partial ^{\alpha }_x\partial
^{\beta }_y ~ \mbox{}_{\sigma } \int_x^{\infty }  \partial ^{\omega
}_z [F(z,y) N(x,z,y)] dz$ \\ for each $|\alpha |+|\beta
|+|\omega |\le m$.
\par In virtue of  Theorems 2.4.1 and 2.5.2 \cite{lucmft12}
or 4.2.5 and 4.2.23 and Corollary 4.2.6 \cite{ludancdnb} there are satisfied the
equalities
$$(8)\quad \sigma _x ~ \mbox{}_{\sigma} \int_x^{\infty }g(z)dz = -g(x)\mbox{  and}$$
$$(9)\quad \mbox{}_{\sigma} \int_{\mbox{}_0x}^x [\sigma _zf(z)]dz
=f(x)-f(\mbox{}_0x)$$ for each continuous function $g$
and a continuously differentiable function $f$, where
$\mbox{}_0x$ is a marked point in $U$,
$$(10)\quad \mbox{}^1\sigma _z ~ \mbox{}_{\sigma } \int_x^{\infty }
F(z,y)N(x,z,y)dz := \sum_{j=0}^{2^r-1} ~ \mbox{}_{\sigma }
\int_x^{\infty } \{ i_j^* [(\partial F(z,y)/\partial z_{\xi
(j)})N(x,z,y)]\psi _j \} dz \mbox{  and}$$
$$(11)\quad \mbox{}^2\sigma _z ~ \mbox{}_{\sigma} \int_{\mbox{}_0x}^x F(z,y)N(x,z,y)dz
:= \sum_{j=0}^{2^r-1} ~ \mbox{}_{\sigma} \int_{\mbox{}_0x}^x \{
i_j^* [ F(z,y)(\partial N(x,z,y)/\partial z_{\xi
(j)})]\psi _j \} dz \mbox{ and}$$
$$(12)\quad ~
\mbox{}^2\sigma _x ~ \mbox{}_{\sigma}\int_x^{\infty }
F (z,y)N(x,z,y)dz := \sum_{j=0}^{2^r-1}
\mbox{}_{\sigma}\int_x^{\infty } \{ i_j^* [F(z,y)(\partial
N(x,z,y)/\partial x_{\xi (j)})] \psi _j \} dz.$$ Therefore, from
Equalities $(8,9)$, 3$(3)$ and 4$(5)$ and Condition $(1)$ we infer that:
$$(13)\quad \sigma _x~ \mbox{}_{\sigma}\int_x^{\infty }  F(z,y) N(x,z,y)dz =
~ \mbox{}^2\sigma _x ~ \mbox{}_{\sigma}\int_x^{\infty }
F (z,y)N(x,z,y)dz - F(x,y)N(x,x,y),$$ since $
F(z,y)N(x,z,y)|_x^{\infty } = - F(x,y)N(x,x,y)$, that
demonstrates Formula $(2)$ for $m=1$ and $A_1=-F(x,y)
N(x,x,y)$. Proceeding by induction for $p=2,...,m$ leads to the identities:
$$(14)\quad \sigma ^p_x~ \mbox{}_{\sigma}\int_x^{\infty }  F(z,y) N(x,z,y)dz
=$$  $$ ~ \sigma _x ~  [\mbox{}^2\sigma ^{p-1}_x
~\mbox{}_{\sigma}\int_x^{\infty } F(z,y)N(x,z,y)dz] +
\sigma _x ~ A_{p-1}(F,N)(x,z,y)$$
$$= ~ \mbox{}^2\sigma ^p_x ~\mbox{}_{\sigma}\int_x^{\infty }
F(z,y)N(x,z,y)dz$$ $$ -  ~ [\mbox{}^2\sigma
^{p-1}_xF(z,y)N(x,z,y)]|_{z=x} + \sigma _x~
A_{p-1}(F,N)(x,y).$$
Thus $(14)$ implies Formulas $(2,4,6)$. Then with the help of
Formulas $(8,9)$ and Condition $(1)$ we infer also that
$$(15)\quad \mbox{}^1\sigma _z~ \mbox{}_{\sigma}\int_x^{\infty }  F(z,y) N(x,z,y)dz =
- ~ \mbox{}^2\sigma _z ~ \mbox{}_{\sigma}\int_x^{\infty }
F(z,y)N(x,z,y)dz +F(z,y)N(x,z,y)|_x^{\infty }$$ $$ =
-F(x,y)N(x,x,y) - ~ \mbox{}^2\sigma _z ~
\mbox{}_{\sigma}\int_x^{\infty } F(z,y)N(x,z,y)dz.$$ Thus
Formulas $(3)$ for $m=1$ and $(7)$ are valid. Then we deduce Formulas
$(3,5)$ by induction on $p=2,...,m$:
$$(16)\quad \mbox{}^1\sigma ^p_z~ \mbox{}_{\sigma}\int_x^{\infty }  F(z,y) N(x,z,y)dz
=$$  $$ \mbox{}^1\sigma ^{p-1}_z [
\mbox{}^1\sigma _z~ \mbox{}_{\sigma}\int_x^{\infty }  F(z,y) N(x,z,y)dz] $$
$$= \mbox{}^1\sigma ^{p-1}_z [ - \mbox{}^2\sigma _z ~
\mbox{}_{\sigma}\int_x^{\infty } F(z,y)N(x,z,y)dz ] -
[\mbox{}^1\sigma ^{p-1}_zF(z,y)N(x,z,y))]|_{z=x}$$
$$= \mbox{}^1\sigma ^{p-2}_z \{ \mbox{}^1\sigma _z [ - \mbox{}^2\sigma _z ~
\mbox{}_{\sigma}\int_x^{\infty } F(z,y)N(x,z,y)dz ] \} $$ $$ -
[\mbox{}^1\sigma ^{p-1}_zF(z,y)N(x,z,y))]|_{z=x}$$
$$= \mbox{}^1\sigma ^{p-2}_z \{ ( - \mbox{}^2\sigma _z)^2 ~
\mbox{}_{\sigma}\int_x^{\infty } F(z,y)N(x,z,y)dz  \} $$
$$+[\mbox{}^1\sigma ^{p-2}_z (  \mbox{}^2\sigma _zF(z,y)N(x,z,y))]|_{z=x}
- [\mbox{}^1\sigma ^{p-1}_zF(z,y)N(x,z,y))]|_{z=x}$$
$$ = ... = (- \mbox{}^2\sigma _z)^p ~
\mbox{}_{\sigma}\int_x^{\infty } F(z,y)N(x,z,y)dz
+ B_p(F(z,y),N(x,z,y))|_{z=x} \mbox{  and}$$
$$B_p(F(z,y),N(x,z,y)) = - (-\mbox{}^2\sigma _z)^{p-1}
F(z,y)N(x,z,y)) + \mbox{}^1\sigma _z B_{p-1}
(F(z,y),N(x,z,y)).$$

\par {\bf 6. Corollary.} {\it If suppositions of Proposition 5 are
satisfied, then
\par $(1)$ $A_2(F,N)(x,y) = - \sigma _x [F(x,y)N(x,x,y)] - ~\mbox{}^2\sigma
_x[F(z,y)N(x,z,y)]|_{z=x}$,
\par $(2)$ $A_3(F,N)(x,y) = - \sigma ^2_x [F(x,y)N(x,x,y)]$\par $ - \sigma _x (~\mbox{}^2\sigma
_x [F(x,y)N(x,z,y)]|_{z=x}) - ~\mbox{}^2\sigma ^2_x[
F(x,y)N(x,z,y)]|_{z=x}$,
\par $(3)$ $A_m(F,N)(x,y) = - \sum_{j=0}^{m-1} ~ \sigma _x ^j \{ [ \mbox{ }^2\sigma _x^{m-1-j} F(z,y)
N(x,z,y) ]|_{z=x} \} $,
\par $(4)$ $B_2(F(z,y),N(x,z,y)) = - \mbox{}^1\sigma _z[F(z,y)N(x,z,y)]
+ ~\mbox{}^2\sigma _z[F(z,y)N(x,z,y)]$,
\par $(5)$ $B_3(F(z,y),N(x,z,y)) = - \mbox{}^1\sigma ^2_z[F(z,y)N(x,z,y)]
$\par  $ + \mbox{}^1\sigma _z(\mbox{}^2\sigma _z[F(z,y)
N(x,z,y)]) - ~\mbox{}^2\sigma ^2_z[F(z,y)
N(x,z,y)]$,
\par $(6)$ $B_m(F(z,y),N(x,z,y)) = [ \sum_{k=0}^{m-1} (-1)^{k+1} \mbox{ }^1\sigma _z^{m-1-k} \mbox{ }^2\sigma _z^{k} ] F(z,y) N(x,z,y)$,
\par $(7)$ $A_2(F,N)(x,y)-B_2(F,N)(x,y) = - 2 ~\mbox{}^2\sigma _x [F(x,y)
N(x,x,y)]$, \\ where $\sigma _xN(x,x,y)=[\sigma _x
N(x,z,y)+\sigma _zN(x,z,y)]|_{z=x}$,
\par $(8)$ $A_3(F,N)(x,y)-B_3(F,N)(x,y)=
- (3 ~\mbox{}^2\sigma ^2_x + ~ \mbox{}^2\sigma _x ~ \mbox{}^2\sigma
_z+ 2~ \mbox{}^2\sigma _z ~ \mbox{}^2\sigma _x) [F(x,y)
N(x,z,y)]|_{z=x}$\par $ - (2 ~\mbox{}^1\sigma _x ~\mbox{}^2\sigma _x
+ ~\mbox{}^2\sigma _x ~\mbox{}^1\sigma _x) [F(x,y)
N(x,x,y)]$.
\par Particularly, if either $p$ is even and $\psi _0=0$, or $F\in Mat_{n\times n}({\bf
R})$ and $N\in Mat_{n\times n}({\cal A}_r)$, then \par $(8)$
$\mbox{}^2\sigma _x^p[F(z,y)N(x,z,y)]= F(z,y)\sigma
_x^pN(x,z,y)$ and $\mbox{}^2\sigma _z^p[F(z,y)
N(x,z,y)]= F(z,y)\sigma _z^pN(x,z,y)$.}
\par {\bf Proof.} From Formulas 5$(13-16)$ Identities $(1-8)$ follow by induction, since
\par $A_m(F,N)(x,y) = - [\mbox{}^2\sigma ^{m-1}_x
F(z,y)N(x,z,y)]|_{z=x}  + \sigma _x A_{m-1}(F,N)(x,y) =$ \par $...=
- [\mbox{}^2\sigma ^{m-1}_x
F(z,y)N(x,z,y)]|_{z=x} - \sigma _x \{ [\mbox{}^2\sigma ^{m-2}_x
F(z,y)N(x,z,y)]|_{z=x} \} $ \par
$-  \sigma ^2_x \{ [\mbox{}^2\sigma ^{m-3}_x
F(z,y)N(x,z,y)]|_{z=x} \} - ... - \sigma ^{m-2}_x \{ [\mbox{}^2\sigma _x
F(z,y)N(x,z,y)]|_{z=x} \} - \sigma ^{m-1}_x
F(x,y)N(x,x,y)$ and
\par $B_m(F(z,y),N(x,z,y)) = - (-\mbox{ }^2 \sigma _z)^{m-1} F(z,y) N(x,z,y) +$\par $ \mbox{ }^1 \sigma _z
B_{m-1}(F(z,y),N(x,z,y)) =...=$
\par $- ~\mbox{}^1\sigma ^{m-1}_z [F(z,y)
N(x,z,y)] + ~\mbox{}^1\sigma ^{m-2}_z(\mbox{}^2\sigma _z[F(z,y)
N(x,z,y)])$ \par $- ~\mbox{}^1\sigma ^{m-3}_z(\mbox{}^2\sigma ^2_z[F(z,y)
N(x,z,y)])+...+(-1)^m (\mbox{}^2\sigma _z)^{m-1} [F(z,y)
N (x,z,y)]$.
\par Particularly when $p$ is even and $\psi _0=0$, $p=2k$, $k\in {\bf
N}$ we get that \par $\sigma
_x^pf(x)=A^kf(x)$ \\ for $p$ times differentiable function $f: U\to
{\cal A}_r$, where \par $Af=\sum_jb_j\partial ^2f(x)/\partial
x_j^2$, $b_j=i_{\xi ^{-1}(j)}^2\in {\bf R}$ according to \S 2.2 \cite{lucmft12}
or Formulas 4.2.4$(7-9)$ \cite{ludancdnb}.
\par On the other hand the
operators $\mbox{}^2\sigma _x^p$ and $\mbox{}^2\sigma _z^p$ commute
with the left multiplication on $F(z,y)\in Mat_{n\times n}({\bf R})$, that is $\mbox{}^2\sigma _x^p[
F(z,y)K(x,z)]= F(z,y)\sigma _x^pK(x,z)$ and
$\mbox{}^2\sigma _z^p[F(z,y)K(x,z)]= F(z,y)\sigma
_z^pK(x,z)$ for $p=2k$, since $\bf R$ is the center of the Cayley-Dickson algebra
${\cal A}_r$.

\section{Some types of integrable nonlinear PDE.}
\par {\bf 1.} Partial differential operators $L_j$ are
considered on domains ${\cal D}(L_j)$ contained in suitable spaces of differentiable functions,
for example, in the space $C^{\infty }(U,Mat_{n\times n}({\cal A}_r))$ of infinitely differentiable by real
variables functions on an open domain $U$ in ${\cal A}_r$ and with values in $Mat_{n\times n}({\cal A}_r)$,
because $U$ has the real shadow $U_{\bf R}$, where $n\in {\bf N}$. Or it is possible to use the Sobolev space
$H^m(U,Mat_{n\times n}({\cal A}_r))$, where $m\ge ord (L_j)$, $ord (L_j)$ denotes the order of a PDE $L_j$,
while on $U$ the Lebesgue measure is provided.
The spaces $C^m(U,Mat_{n\times n}({\cal A}_r))$ and $H^m(U,Mat_{n\times n}({\cal A}_r))$ with $m\le \infty $ are linear
over the real field $\bf R$, also they have the structure of
the left and the right modules over the Cayley-Dickson algebra
${\cal A}_r$, $r\ge 2$.
To each Cayley-Dickson number $z = z_0i_0+...+z_{2^r-1}i_{2^r-1} \in {\cal A}_r$ there corresponds a vector $[z]=(z_0,...,z_{2^r-1})$ in its real shadow ${\bf R}^{2^r}$, where $z_j\in \bf R$ for each $j$.
For functions $f([z])$ of $[z]$ we shall write for short $f(z)$ also.
\par Henceforth, if something other will not be specified, we shall take
a function $N$ may be depending on $F$, $K$ and satisfying the following conditions:
\par $(1)$ $N(x,y) = EK(x,y)$ with an operator $E$ in the form
\par $(2)$ $E= BST_g$,
\par $(3)$ $[L_j, E] =0$ for each $j$,
\\ where $B$ is a nonzero bounded right ${\cal A}_r$ linear (or strongly right ${\cal A}_r$
linear) operator,
$~ S=S(x,y)\in Aut (Mat_{n\times n}({\cal A}_r))$, so that $B$ is independent of $x, y\in U$,
$g\in Diff^{\infty }(U^2_{\bf R})$, $g=(g_1,g_2)$, $g_l(U^2_{\bf R})=U_{\bf R}$ for
$l=1$ and $l=2$, $U_{\bf R}$ denotes the real shadow of the domain $U$,
$Aut (Mat_{n\times n}({\cal A}_r))$ notates the automorphism group of the algebra $Mat_{n\times n}({\cal A}_r)$,
\par $(4)$ $T_gK(x,y) := K(g_1(x,y),g_2(x,y))$.
\par Condition $(3)$ is implied by the following:
\par  $(5)$ $[L_j, B] =0$ for each $j$,
\par $(6)$ $L_{j,x,y} (T_gK(x,y)) = T_g (L_{j,x,y}K(x,y))$ and
\par $(7)$ $L_{j,x,y} (S(x,y)K(x,y)) = S(x,y)(L_{j,x,y}K(x,y))$ for each $j$
and each $x, y \in U$, where $L_{j,x,y}$ are PDOs considered below.
\par Evidently Conditions $(6,7)$ are fulfilled,
when $L_{j,x,y}$ are polynomials of $\sigma _x^k$ and $\sigma _y^k$, all coefficients of $L_j$ are real and
the following stronger conditions are imposed:
\par $(8)$ $\sigma _x^k(T_gK(x,y)) = T_g (\sigma _x^kK(x,y))$, $\sigma _y^k(T_gK(x,y)) = T_g (\sigma _y^kK(x,y))$
and
\par $(9)$ $\sigma _x^k(S(x,y)K(x,y)) = S(x,y)(\sigma _x^kK(x,y))$, $\sigma _y^k(S(x,y)K(x,y)) = S(x,y)(\sigma _y^kK(x,y))$, \\ since $S|_{i_0{\bf R}}=I$.
\\ If coefficients of $L_j$ may be Cayley-Dickson numbers, $S=I$, then Condition $(8)$
will suffice as well.
\par Particularly there may be $E=B$, $B\in SL_n({\bf R})$, or $E=I$. It will also be indicated,
when $E$ or $K$ and hence $N$ depend on some parameter or a variable.

\par {\bf 2. General approach to solutions of nonlinear vector partial
differential equations with the help of non-commutative integration
over Cayley-Dickson algebras.} We consider an equation over the
Cayley-Dickson algebra ${\cal A}_r$ which is presented in the non-commutative
line integral form:
$$(1)\quad K(x,y) = F(x,y) + p ~ \mbox{}_{\sigma }
\int_x^{\infty } F(z,y)N(x,z,y)dz ,$$ where $K$, $
F$ and $N$ are continuous integrable functions of ${\cal A}_r$
variables $x, y, z \in U$ so that $F$, $K$ and $N$ have values in $Mat_{n\times n}({\cal A}_r)$,
where $n\ge 1$, $r\ge 2$, $ ~ N$ and $K$ are related by 1$(1,2)$,
$p ~ \in {\bf R}\setminus \{ 0 \} $ is a non-zero real
constant. These functions $F$, $K$ and $N$ may depend on additional parameters $t,
\tau,...$. \par  At first it is necessary to specify the
function $N$ and its expression throughout $F$ and $
K$. It is supposed that an operator
\par $(2)$ $(I- {\sf A}_xE)K(x,y)=F(x,y)$ is invertible, \\ when
$N(x,z,y) = E_yK(x,z)$ for each $x, y, z \in U$,
so that $(I- {\sf A}_xE)^{-1}$ is continuous, where $I$ denotes the
unit operator, \par $(3)$ $ {\sf A}_x K(x,y) :=  p ~ \mbox{}_{\sigma }
\int_x^{\infty } F(z,y)K(x,z)dz $ \\ is an operator
acting by variables $x$. \par Then ${\bf
R}$-linear partial differential operators $L_k$ over
the Cayley-Dickson algebra ${\cal A}_r$ are provided for $k=1,...,k_0$, where $k_0\in {\bf N}$.
It is frequently helpful to consider their decompositions:
\par $(4)$ $L_kf = \sum_j i_j^* (L_{k,j}f)$, \\ where
$f$ is a differentiable function in the domain of each operator $L_k$,
$~L_{k,j}$ are components of the operators $L_k$ so that each $L_{k,j}$
is a PDO written in real variables with real coefficients. That is $L_{k,j}g$
is a real-valued function for each $ord (L_{k,j})$ times differentiable
real-valued function $g$ in the domain of $L_{k,j}$ for every $j$,
where $ord (L_{k,j})$ denotes the order of the PDO $L_{k,j}$.
Next the conditions are imposed on the function $F$:
\par $(5)$ $L_kF =0$ \\ for  $k=1,...,k_0$. \par It may be necessary to consider in some
problems stronger conditions:
\par $(6)$ $\sum_{j\in \Psi _l} i_j^* [~ c_{k,j} (L_{k,0}F)+L_{k,j}F] =0$ \\
for each $k$ and $1\le l \le m$, where $c_{k,j}$
are constants $c_{k,j}\in {\cal A}_r$, $~\Psi _l\subset \{
0,1,...,2^r-1 \} $ for each $l$, $ ~ \bigcup_l \Psi _l = \{ 0, 1,
..., 2^r-1 \} $, $~ \Psi _n \cap \Psi _l = \emptyset $ for each
$n\ne l$, $~1\le m\le 2^r$.  There is not excluded that the
coefficients $c_{k,j}$ or the operators
$L_{k,j}$ may be zero for some $(k,j)$. \par After this a
function $K$ is determined from Equation $(2)$. \par This function $K$ may
be satisfying some PDEs, when suitable PDOs $L_s$ and the operator $E$ are chosen
(see also \S 1). Indeed acting by
the operator $L_k$ from the left on both sides of $(2)$ one may get
with the help of Conditions
either $(5)$ or $(6)$ the PDEs
either
\par $(7)$ $L_s[(I-{\sf A}_xE_y)K]=0$ or
\par $(8)$ $\sum_{j\in \Psi _k} i_j^*
\{ ~ c_{k,j} L_{s,0}[(I- {\sf A}_xE_y)K] + L_{s,j}[(I-
{\sf A}_xE_y)K]\} =0$ for each $k=1,...,m$ respectively for $s=1,...,k_1$, where $k_1\le k_0$.
\\ Therefore this leads to the equalities
\par $(9)$ $(I- {\sf A}_xE_y)(L_sK)=R_s(K)$ for $s=1,...,k_1$, \\ where each operator of the form
\par $(10)$ $R_s(f)=(I- {\sf A}_xE_y)(L_sf)-L_s[(I- {\sf A}_xE_y)f]$ \\
is obtained by calculations of appearing commutators $[A,B]=AB-BA$ and
anti-commutators $\{ A,B \} = AB+BA$ of operators $(I-(
{\sf A}_xE_y)_0)$, $~ ({\sf A}_xE_y)_j$, $~L_{s,j}$, $~ j=0,...,2^r-1$.
The latter can be realized when the function $N$ and the PDOs
$L_j$ are chosen such that
\par $(11)$ $R_s(K)=(I- {\sf A}_xE_y)M_s(K)$ for $s=1,...,k_1$, \\ where $M_s(K)$ are operators
or functionals acting on $K$. Generally the operators $M_s$ may be non-${\bf R}$-linear
and besides terms of a partial differential operator it may contain terms containing
the integral operator ${\sf A}$.
Therefore due to Condition $(2)$ the function $K$ must satisfy the PDEs or the partial integro-differential
equations (PIDEs)
\par $(12)$ $L_sK-M_s(K)=0$ for $s= 1,...,k_1$, \\ which generally may be non-${\bf
R}$-linear. Thus each solution $K$ of the $\bf R$-linear
integral equation $(1)$ should also be the solution of the aforementioned
PDEs or PIDEs $(12)$. \par It is worthwhile to choose the ${\cal A}_r$ vector
independent PDOs $L_{s}$ for $s=1,...,k_0$ and so that $c_{k,j}\in {\bf R}i_{\xi (k,j)}$
and $\xi (k,j)\in \{ 0,1,..,2^r-1 \} $ for each $k, j$.
\par Henceforward, if something other will not be outlined, we consider the variants:
\par $(13)$ $F, K, N \in Mat_{n\times n}({\cal A}_r)$ with $2\le r \le 3$ and $B$ is the strongly right ${\cal A}_r$-linear operator; or
\par $(14)$ $F\in Mat_{n\times n}({\bf R})$ and $K, N \in Mat_{n\times n}({\cal A}_r)$ with $2\le r$
and $B$ is the right ${\cal A}_r$-linear operator (see also \S 1), where $1\le n\in {\bf N}$.

\par {\bf 3. Theorem.} {\it Suppose that conditions of Proposition 2.5 and 2.2$(RS)$ are fulfilled
over the Cayley-Dickson algebra ${\cal A}_r$ with $2\le r$ and on a domain
$U$ satisfying Conditions 2.1$(D1,D2)$
for the corresponding terms of operators $L_s$ for all $s=1,...,k_0$
so that
\par $(1)$ the appearing in the terms $M_s(K)$
integrals uniformly converge by parameters on compact sub-domains in $U$
and
\par $(2)$ $\lim_{z\to \infty } \partial_x^{\alpha }\partial_y^{\beta } \partial_z^{\omega }
(F(z,y)N(x,z,y)]=0$ \\ the limit converges uniformly by $x, y \in U\setminus V$
for some compact subset $V$ in $U$ and for each $|\alpha |+|\beta |+|\omega |\le m$,
where $1\le m=\max \{ deg (L_s): s=1,...,k_0 \} $ and
\par $(3)$ the operator $(I-{\sf A}_xE_y)$ is invertible,
where $F$ is in the domain of PDOs $L_1,...,L_{k_0}$,
$ ~ F(x,y)\in Mat_{n\times n}({\cal A}_r)$ and $K(x,y)\in Mat_{n\times n}({\cal A}_r)$,
$n\in {\bf N}$.
\par Then there exists a solution $K$ of PDEs or PIDEs
2$(12)$ such that $K$ is given by Formulas 1$(1,2)$, 2$(1)$ and either 2$(5)$ or 2$(6)$.}
\par {\bf Proof.} The anti-derivative operator $g\mapsto ~\mbox{}_{\sigma }
\int_{\mbox{}_0x}^x g(z)dz$ is compact from $C^0(V,{\cal A}_r)$ into
$C^0(V,{\cal A}_r)$ for a compact domain $V$ in ${\cal A}_r$, where
$C^0(V,{\cal A}_r)$ is the Banach space over ${\cal A}_r$ of all
continuous functions $g: V\to {\cal A}_r$ supplied with the supremum
norm $\| g \| := \sup_{x\in V} |g(x)|$, $~\mbox{}_0x$ is a marked
point in $V$, $ ~ x\in V$. A function $F$ satisfying
the system of ${\bf R}$ linear PDEs 2$(5)$
or 2$(6)$ is continuous. \par Therefore, due to Conditions $(1-3)$ the anti-derivative operator $\mbox{}_{\sigma }
\int_x^{\infty } F(z,y)N(x,z,y)dz $ is compact. Hence there exists $\delta >0$
such  that the operator $I-{\sf A}_xE_y$ is invertible when $|p|<\delta $, where $p\in {\bf R}\setminus \{ 0 \} $.
Mention that the operator $T_g$ is strongly left and right ${\cal A}_r$-linear (see \S 1), while
$S$ is the automorphism of the Cayley-Dickson algebra, that is
$S[ab] = S[a]S[b]$ and $S[a+b] = S[a]+S[b]$ for each $a, b \in {\cal A}_r$.
\par Since the operator $(I-{\sf A}_xE_y)$ is invertible and Conditions 2.2$(RS)$
and either 2$(13)$ or 2$(14)$ are satisfied, then Equation 2$(12)$ can be resolved:
\par $(4)$ $\sum_k ~\mbox{}_kf(x,y)~\mbox{}_kg(y)=(I-{\sf A}_xE_y)^{-1}u(x,y)$, \\ since
if $A: X\to X$ is a bounded ${\bf R}$ linear operator on a
Banach space $X$ with the norm $ \| A \| <1$, then the inverse of $I-A$ exists:
\par $(I-A)^{-1} = \sum_{n=0}^{\infty } A^n$.
Applying Proposition 2.5 and \S 2 we get the statement of this theorem.

\par {\bf 4. Remark.} If Condition 2.2$(RS)$ is not fulfilled,
the corresponding system of PDEs in real components $({\sf A}_xE_y)_{j,s}$, $~\mbox{}_kf_s$ and
$~\mbox{}_kg_s$ can be considered.

\par {\bf 5. Lemma.} {\it Let suppositions of Proposition 2.5 be satisfied
and the operator ${\sf A}_x$ be given by Formula 2$(3)$, let also
$E=E_y$ may be depending on the parameter $y\in U$ and let $N(x,z,y) = E_y K(x,z)$
(see Formulas 1$(1-4)$). Suppose that $F(x,y)\in Mat_{n\times n}({\bf R})$
and $K(x,z)\in Mat_{n\times n}({\cal A}_r)$ for each $x, y, z \in U$, where $r\ge 2$.
Then
\par $(1)$ $A_m(F, E_yK)(x,y) = (I-{\sf A}_xE_y) \hat{A}_m (K, E_yK)(x,y) + P_m (K, E_yK)(x,y)$,
\par $(2)$ $B_m(F, E_yK)(x,y) = (I-{\sf A}_xE_y) \hat{B}_m (K, E_yK)(x,y) + Q_m (K, E_yK)(x,y)$, where
$$(3)\quad \hat{A}_{m,x,y} (K(z,y), E_yK(x,z))|_{z=x} = \hat{A}_m (K, E_yK)(x,y) = - \sum_{j=0}^{m-1} \sigma _x^j K_{1,m-j-1}(x,y)$$
$$ + p\sum_{j=1}^{m-1} \sum_{j_1=0}^{j-1}\sigma _x^{j_1}  K_{2,m-j-1,j-j_1-1}(x,y)$$ $$ + p^2\sum_{j=1}^{m-1} \sum_{j_1=1}^{j-1}\sum_{j_2=0}^{j_1-1}\sigma _x^{j_2}  K_{3,m-j-1,j-j_1-1,j_1-j_2-1}(x,y) +...+
p^{m-2} K_{m-1,0,...,0}(x,y),$$
$$(4)\quad \hat{B}_m (K(z,y), E_yK(x,z)) = \sum_{j=1}^{m} (-1)^j \{ \mbox{ }^1\sigma _z^j [K(z,y)
\sigma _v^{j-1} (E_yK(w,v))]$$ $$ + p \hat{A}_{m-j,z,y} (K(z,y), E_yK(z,x) (\sigma _v^{j-1}
E_yK(w,v)))\} |_{v=z, w=x} ,$$
\par $(5)$ $\hat{B}_m (K, E_yK)(x,y) := \hat{B}_m (K(z,y), E_yK(x,z))|_{z=x}$,
\par $(6)$ $\hat{A}_1 (K, E_yK)(x,y) = - K(z,y) (E_yK(x,z))|_{z=x}$,
\par $(7)$ $\hat{B}_1 (K(z,y), E_yK(x,z)) = - K(z,y) (E_yK(x,z))$
\\ for each $m\ge 2$ in $(3,4)$, where $\sum_{j=l}^m a_j:=0$ for all $l>m$,
\par $(8)$ $K_{1,j} (x,y|z) := K(x,y) \sigma_x^{j} (E_yK(x,z))$,
\par $(9)$ $K_{m,l_1,...,l_m} (x,y|z) := K(x,y) [\sigma_x^{l_m} E_yK_{m-1,l_1,...,l_{m-1}} (x,z)]$,
\par $(10)$ $ K_{m,l_1,...,l_m} (x,y) := K_{m,l_1,...,l_m} (x,y|z)|_{z=x}$,
\par $(11)$ $P_m(K(z,y), E_yK(x,z))|_{z=x} = P_m(K, E_yK)(x,y) :=
{\sf A}_x \{ \sum_{j=1}^{m-1} [\sigma_x^j, E_y] K_{1,m-j-1} (x,y)$ \par $ + p \sum_{j=1}^{m-1} \sum_{j_1=1}^{j-1} [\sigma_x^{j_1}, E_y]
K_{2,m-j-1,j-j_1-1} (x,y) +...$ \par $+ p^{m-3} \sum_{j=1}^{m-1} \sum_{j_1=1}^{j-1} \sum_{j_{m-3}=1}^{j_{m-4}-1} [\sigma_x^{j_{m-3}}, E_y]
K_{m-2,m-j-1,j-j_1-1,...,j_{m-4}-j_{m-3}-1} (x,y) \} $,
\par $(12)$ $Q_m(K,EK)(x,y|z) := $ \par $\sum_{j=1}^{m-1} (-1)^j P_{m-j}(K(\eta ,y), E_yK(z,\eta )(\sigma _v^{j-1}E_yK(w,v))|_{\eta =z, v=z, w=x}$,
\par $(13)$ $Q_m(K,EK)(x,y) := Q_m(K,EK)(x,y|z)|_{z=x}$. }
\par {\bf Proof.}  Formulas $(6)$ and $(7)$ follow immediately from that of
2.2$(2,3)$ and 2.5$(6,7)$.
Write $A_m$ for each $m\ge 2$ in the form:
$$(14)\quad A_m (F,E_yK)(x,y) = - \sum_{j=1}^m \sigma_x^{j-1} \{ [ \mbox{ }^2 \sigma_x^{m-j}F(z,y)
(E_yK(x,z))]|_{z=x} \} ,$$ where $\sigma ^0=I$. Using that $F(z,y) = (I-{\sf A}_zE_y)K(z,y)$ and
$F(z,y)\in Mat_{n\times n}({\bf R})$ we get from $(14)$:
$$(15)\quad A_m (F,E_yK)(x,y) = - \sum_{j=1}^m \sigma_x^{j-1} \{ [(I-{\sf A}_xE_y)K(x,y)
(\sigma_x^{m-j}E_yK(x,z))]|_{z=x} \}.$$
In virtue of Proposition 2.5 we deduce from $(12)$ that
$$(16)\quad A_m(F,E_yK)(x,y) = -  (I-{\sf A}_xE_y) \{  \sum_{j=1}^m \sigma_x^{j-1} K_{1,m-j}(x,y) \}  +$$
$$p \sum_{j=2}^mA_{j-1} (F(z,y),E_yK_{1,m-j}(x,z))|_{z=x} + {\sf A}_x \sum_{j=1}^{m-1} [\sigma _x^j, E_y] K_{1,m-j-1}(x,y) $$
$$ = -  (I-{\sf A}_xE_y) \{  \sum_{j=0}^{m-1} \sigma_x^{j} K_{1,m-j-1}(x,y)  + p \sum_{j=1}^{m-1} \sum_{j_1=0}^{j-1}\sigma _x^{j_1}  K_{2,m-j-1,j-j_1-1}(x,y) \} $$
$$+ {\sf A}_x \sum_{j=1}^{m-1} [\sigma _x^j, E_y] K_{1,m-j-1}(x,y) + p {\sf A}_x \sum_{j=1}^{m-1} \sum_{j_1=1}^{j-1} [\sigma _x^{j_1}, E_y] K_{2,m-j-1, j- j_1 -1}(x,y)$$
$$ + p^2  \{  \sum_{j=1}^{m-1} \sum_{j_1=1}^{j-1} A_{j_1} (F(z,y),E_yK_{2,m-j-1,j-j_1-1}(x,z))|_{z=x} =... ,$$
since $[\sigma _x^j,E_y] + E_y \sigma _x^j = \sigma _x^j E_y$ for $j\ge 1$,
\par $A_xE_yK(x,y) = p \mbox{ }_{\sigma } \int_x^{\infty } F(z,y)E_yK(x,z)dz$.
Iterating relations $(13)$ we infer by induction Formulas $(1,3,11)$.
Then we have
$$(17)\quad B_m(F(z,y),E_yK(x,z)) = \sum_{j=1}^m (-1)^j \mbox{ }^1\sigma _z^{m-j} \mbox{ }^2\sigma _z^{j-1}
F(z,y) (E_yK(x,z))$$ and for $F(z,y)\in Mat_{n\times n}({\bf R})$ for each $z, y\in U$
this reduces to:
$$(18)\quad B_m(F(z,y),E_yK(x,z)) = \sum_{j=1}^m (-1)^j \mbox{ }^1\sigma _z^{m-j}
F(z,y) (\sigma _z^{j-1}(E_yK(x,z))) $$
$$= \sum_{j=1}^m (-1)^{j_1} \sigma _z^{m-j} (I-{\sf A}_z E_y)K(z,y) (\sigma _v^{j-1} (E_yK(w,v))|_{v=z, w=x}).$$
Therefore, in view of Proposition 2.5 and Formula $(1)$ the identity
$$(19)\quad B_m(F(z,y),E_yK(x,z)) = (I-A_zE_y) \{ \sum_{j=1}^m (-1)^j \{ \mbox{ }^1\sigma _z^{m-j}
K(z,y) (\sigma _v^{j-1}(E_yK(w,v)))$$ $$ + p \hat{A}_{m-j,z,y} (K(\eta ,y), E_yK(z,\eta ) (\sigma _v^{j-1}
E_yK(w,v))) \} \} |_{\eta =z, v=z, w=x} $$
$$+ \sum_{j=1}^{m-1} (-1)^j P_{m-j} (K(\eta ,y), E_y K(z,\eta )(\sigma _v^{j-1} E_yK(w,v)) |_{\eta =z, v=z, w=x}$$
is valid, since $$A_zE_yK(z,y) \{ \sigma _v^{j-1}(E_yK(w,v)) \} |_{v=z, w=x} = $$   $$
p \mbox{ }_{\sigma } \int_z^{\infty } F(\eta ,y) \{ E_yK(z,\eta ) \{ \sigma _v^{j-1}(E_yK(w,v)) \} \} d\eta |_{v=z, w=x} .$$

\par {\bf 6. Proposition.} {\it Suppose that \par $(1)$ a PDO $L_j$ is a polynomial $\Omega _j(\sigma _x, \sigma _y)$
of $\sigma _x$ and $\sigma _y$ for each $j=1,...,k_0$, coefficients of $\Omega _j$ are real and
Condition 1$(3)$ is fulfilled for all $j$;
\par $(2)$ $L_{s,x,y}K(x,y) - ({\sf A}_x E_y L_{s,x,y})K(x,y) =: R_s(K)(x,y)$ for each $s=1,...,k_1$;
\\ $F(x,y)$ is in $Mat_{n\times n}({\bf R})$ and $K(x,y) \in Mat_{n\times n}({\cal A}_r)$
for each $x, y \in U$ (see 2$(3)$), $\sigma $ and ${\sf A}$ are
over the Cayley-Dickson algebra ${\cal A}_r$, $2\le r$, $n\in {\bf N}$, $1\le k_1\le k_0$;
\par $(3)$ $L_{j,x,y} F(x,y)=0$ for every $x, y \in U$ and $j=1,...,k_0$.
\par Then there exists a polynomial $M_s$ of $K, E, \sigma $ and ${\sf A}$ such that
\par $(4)$ $R_s(K)(x,y) = (I-{\sf A}_x E_y) M_s(K)(x,y)$ for each $s=1,...,k_1$.}
\par {\bf Proof.} Proposition 2.5 and Corollary 2.6 imply that $R_s(K)(x,y)$ can be expressed as
a polynomial of $A_m(F,EK)$, $B_m(F,EK)$, $\sigma $ and ${\sf A}EK$,
where $m\in \bf N$.
\par Take an algebra ${\cal B}$ over the real field
generated by the operators $\sigma $, ${\sf A}$, $E$ and $I$:
\par ${\cal B} = alg_{\bf R} (\sigma _x, \sigma _y, \sigma _z, {\sf A}_x, {\sf A}_y, {\sf A}_z,
E_y(x,z) ~ \forall  x, y, z\in U; I)$, where $I$ denotes the unit operator.
In view of Proposition 3.1 \cite{lucveleq2013}, Theorems 2.4.1 and 2.5.2 \cite{lucmft12}
the algebra ${\cal B}$ is associative, since $F(z,y)$ is in
$Mat_{n\times n}({\bf R})$ for each $z, y \in U$, the algebra $Mat_{n\times n}({\bf R})$ is associative,
also $E$ is given by 1$(2,4)$. Therefore, there exists the Lie algebra ${\cal L} ({\cal B})$
generated from ${\cal B}$ with the help of commutators $[H,G]:= HG-GH$ of elements $H, G \in
{\cal B}$ (see also about abstract algebras of operators and their Lie algebras in
\cite{schaferb}). Then $Y_s := [L_s, E]{\cal L}({\cal B})$ is the (two-sided) ideal in ${\cal L}({\cal B})$
and hence there exist the quotient algebra ${\cal L}_s := {\cal L}({\cal B})/Y_s$
and the quotient morphism $\pi _s: {\cal L}({\cal B})\to {\cal L}_s$.
\par Next consider the universal enveloping algebra ${\cal U}$ of the Lie algebra
${\cal L}({\cal B})$. In virtue of Proposition 2.1.1 \cite{bourbalglie} there exists a unique homomorphism
$\tau $ from ${\cal U}$ into ${\cal B}$. The algebra $C^{\infty }(U,Mat_{n\times n}({\cal A}_r))$ over the real field
also has the structure of the left module
of the operator ring  ${\cal B}$ and hence of ${\cal L}({\cal B})$ and ${\cal U}$ as well, where
$C^{\infty }(U,Mat_{n\times n}({\cal A}_r))$ denotes the algebra of all infinitely differentiable functions
from $U_{\bf R}$ into $Mat_{n\times n}({\cal A}_r)$ (see \S 1).
Since $C^{\infty }(U,Mat_{n\times n}({\cal A}_r))$ is dense in $C^l(U,Mat_{n\times n}({\cal A}_r))$, then it is
sufficient to consider $C^{\infty }(U,Mat_{n\times n}({\cal A}_r))$, where $l=\max_{s=1,...,k_0} ord (L_s)$.
On the other hand, $[L_s, E]{\cal U} =: {\cal U}_s$ is the (two-sided) ideal in
${\cal U}$. \par Let
${\cal P}(x,y)$ denote the $\bf R$-linear algebra generated by sums and products of all terms $QP$ so that
$P$ are polynomials of functions
$K\in C^{\infty }(U,Mat_{n\times n}({\cal A}_r))$ and $Q$ are acting on them polynomials of operators
$\sigma $, $\sf A$, $E$ (or $B$, $S$, $T_g$ instead of $E$, since $E=BST_g$), where coefficients of $P$ and $Q$
are chosen to be real, since coefficients of each polynomial $\Omega _j$ are real.
Certainly the equality $\alpha T=T\alpha $ is valid for each $T\in Mat_{n\times n}({\bf R})$ and $\alpha \in {\cal A}_r$,
since $\bf R$ is the center of the Cayley-Dickson algebra ${\cal A}_r$, $2\le r$
and $(\alpha T)_{i,j} = \alpha T_{i,j}= T_{i,j}\alpha = (T\alpha )_{i,j}$ for each $(i,j)$ matrix element $(\alpha T)_{i,j}$ of
$\alpha T$.
\par The polynomial $R_s(K)(x,y)$ is calculated with the help of Conditions
$(3)$, where $s=1,...,k_1$ (see also \S 2). From Formula $(2)$ it follows that the polynomial $R_s(K)(x,y)$ belongs to
$(I-{\sf A}_x E_y) {\cal P}(x,y) + \sum_{j=1}^{k_0} \{ {\cal U}_j C^{\infty }(U,Mat_{n\times n}({\cal A}_r)) \} $ for each
$s=1,...,k_0$.
Applying the quotient mapping $\pi _j$ for all $j=1,...,k_0$ and using Proposition 2.3.3 \cite{bourbalglie} we get Formula $(4)$,
since $\pi _j ({\cal U}_j) = 0$ and Condition 1$(3)$ is imposed for each $j=1,...,k_0$.

\par {\bf 7. Example.} Take two partial differential operators
\par $(1)$ $L_1=L_{1,x,y} := \sum_l a_l ((-\sigma _x)^l - \sigma _y^l)$ and
\par $(2)$ $L_2=L_{2,x,y} := \sum_l a_l (\sigma _x^l - ( - \sigma _y)^l)$,
\\ where $a_l\in {\bf R}$ for each $l$ when $\sigma $ is over ${\cal A}_r$
with $r\ge 2$, the sum is finite or infinite, $l\in {\bf N}$.
The functions $F(x,y)$ and $K(x,y)$ of ${\cal A}_r$
variables $x, y\in U$ have values in $Mat_{n\times n}({\bf R})$ and $Mat_{n\times n}({\cal A}_r)$
respectively, where $n\ge 1$, $r\ge 2$. A domain $U$
in ${\cal A}_r$ satisfies conditions 2.1$(D1,D2)$ with $\infty
\in U$.
On a function $F(x,y)$ are imposed two conditions:
\par $(3)$ $L_{1,x,y}F(x,y)=0$ and
\par $(4)$ $L_{2,x,y}F(x,y)=0$.
\par Suppose that conditions of Proposition 2.5 are fulfilled and
\par $(5)$ $K(x,y) = F(x,y) + p \mbox{}_{\sigma} \int_x^{\infty }F(z,y)N(x,z)dz$,
\\ where $p$ is a non-zero real parameter, $N(x,z) = E K(x,z)$ for each $x, z\in U$,
while $E$ is a bounded right ${\cal A}_r$-linear operator satisfying Conditions 1$(2-4)$ and either 2$(13)$
or 2$(14)$.
\par Condition $(3)$ is equivalent to
\par $(6)$ $\sum_l a_l (-\sigma _x)^lF(x,y) = \sum_l a_l \sigma _y^l F(x,y)$ and $(4)$ to
\par $(7)$ $\sum_l a_l \sigma _x^lF(x,y) = \sum_l a_l (-1)^l\sigma _y^l F(x,y)=0$ correspondingly.
Acting on both sides of the equality $(5)$ by the operator $L_1$ and using $(6)$
and Proposition 2.5 we get
\par $(8)$ $L_{1,x,y} K(x,y) = p \sum_l a_l ((-\sigma _x)^l - (-\mbox{ }^1\sigma _z)^l) \mbox{}_{\sigma} \int_x^{\infty }F(z,y)N(x,z)dz$   \par $= p \sum_l a_l ((-\mbox{ }^2\sigma _x)^l -  \mbox{ }^2\sigma _z^l) \mbox{}_{\sigma} \int_x^{\infty }F(z,y)N(x,z)dz$
\par $ + p \sum_l a_l (-1)^l (A_l(F;N)(x,y) - B_l(F;N)(x,y))$ and hence
\par $(9)$ $L_{1,x,y} K(x,y) = p \mbox{ }^2L_{1,x,z} \mbox{}_{\sigma} \int_x^{\infty }F(z,y)N(x,z)dz$
\par $ + p \sum_l (-1)^l a_l (A_l(F;N)(x,y) - B_l(F;N)(x,y))$.
\par Then from $(5,7)$ we infer that
\par $(10)$ $L_{2,x,y} K(x,y) = p \sum_l a_l (\sigma _x^l - \mbox{ }^1\sigma _z^l) \mbox{}_{\sigma} \int_x^{\infty }F(z,y)N(x,z)dz$   \par $= p \sum_l a_l (\mbox{ }^2\sigma _x^l - (-\mbox{ }^2\sigma _z)^l) \mbox{}_{\sigma} \int_x^{\infty }F(z,y)N(x,z)dz$
\par $ + p \sum_l a_l (A_l(F;N)(x,y) - B_l(F;N)(x,y))$ and consequently,
\par $(11)$ $L_{2,x,y} K(x,y) = p \mbox{ }^2L_{2,x,z} \mbox{}_{\sigma} \int_x^{\infty }F(z,y)N(x,z)dz$
\par $ + p \sum_l a_l (A_l(F;N)(x,y) - B_l(F;N)(x,y)).$
\par Then Equalities $(9,11)$ imply that
\par $(12)$ $ (L_{1,x,y} \pm  L_{2,x,y} ) K(x,y) = p (\mbox{ }^2L_{1,x,z} \pm  \mbox{ }^2L_{2,x,z})\mbox{}_{\sigma} \int_x^{\infty }F(z,y)N(x,z)dz$
\par $ + p \sum_l a_l (\pm 1 +(-1)^l) (A_l(F;N)(x,y) - B_l(F;N)(x,y)).$
\par We take into account sufficiently small values of the parameter $p$, when the operator $I-{\sf A}_xE$ is invertible,
for example, $ \| {\sf A}_xE \| <1$, where
\par ${\sf A}_xK(x,y)  = p \mbox{}_{\sigma} \int_x^{\infty }F(z,y)K(x,z)dz$.
In the case $F\in Mat_{n\times n} ({\bf R})$ and $K\in Mat_{n\times n} ({\cal A}_r)$
with $\sigma $ over
${\cal A}_r$, from $(5,12)$, Lemma 5 and Proposition 6 it follows that $K$ satisfies
the nonlinear PDE
\par $(13)$ $L_{x,y}^{\pm } K(x,y) - p \sum_l (\pm 1+(-1)^l)a_l [\hat{A}_l(K;EK)(x,y) -  B_l({\hat K};EK)(x,y)]=0,$ where
\par $(14)$ $L_{x,y}^{\pm }  = \sum_l (\pm 1+(-1)^l)a_l (\sigma _x^l - \sigma _y^l),$
\\ $\hat{A}_l$ and $\hat{B}_l$ are given by Formulas 5$(3-7)$, since
$<a,b,c> =0$ when particularly $a\in Mat_{n\times n}({\bf R})$, where $<e,b,c> = (eb)c - e(bc)$
denotes the associator of the Cayley-Dickson matrices $e, b, c \in Mat_{n\times n}({\cal A}_r)$, also since
$\alpha a= a\alpha $ for each $\alpha \in {\cal A}_r$.
 A solution of $(13)$ reduces to linear PDEs
and is prescribed by $(3-5)$. Equivalently the function $K$ satisfies also the PDEs
\par $(15)$ $L_{s,x,y} K(x,y) - p \sum_l (-1)^{sl}a_l [\hat{A}_l(K;EK)(x,y) -  \hat{B}_l(K;EK)(x,y)]=0$ \\
for $s=1$ and $s=2$. Instead of this system it is possible also to consider
separately PDOs $L_1$ and $L_2$ and the corresponding PDEs for $F$ and $K$
as well. Thus with the help of Theorem 3 we get the following.

\par {\bf 7.1. Theorem.} {\it Suppose that conditions of Theorem 3 and Example 7 are fulfilled,
then a solution of PDE $(15)$ is given by $(3,5)$, where $L_1$ is prescribed by Formula $(1)$, $s=1$.}

\par {\bf 8. Example.} Let PDOs be
\par $(1)$ $L_1 =L_{1;x,y}= \sigma _x - \sigma _y$,
\par $(2)$ $L_{2,j} =L_{2,j;x,y}= \sum_l (a_l \sigma _x^l +(-1)^{jl} b_l \sigma _y^l)$,
\\ where $a_l, b_l \in {\bf R}$ for each $l$ when $\sigma $ is over ${\cal A}_m$
with $m\ge 2$, the sum is finite or infinite, $j=1$ or $j=2$. It is also supposed that the functions
$F(x,y)$ and $K(x,y)$ of ${\cal A}_m$
variables $x, y\in U$ have values in $Mat_{n\times n}({\bf R})$ and $Mat_{n\times n}({\cal A}_m)$
correspondingly, where $n\ge 1$, $m\ge 2$. A domain $U$
in ${\cal A}_m$ satisfies Conditions 2.1$(D1,D2)$ with $\infty
\in U$. Suppose that \par $(3)$ $L_{1;x,y}F(x,y) =0$ and consider the integral relation:
\par $(4)$ $K(x,y) = F(x,y) + p \mbox{}_{\sigma} \int_x^{\infty }F(z,y)N(x,z)dz$, \\ where
$N$ and $K$ are related by Formulas 1$(1,2)$.
\par Using condition $(3)$ we can write $F(x,y) = F(\frac{x+y}{2})$.
Therefore we deduce that
$$(5)\quad L_{2,j;x,y}K(x,y) = L_{2,j;x,y}F(\frac{x+y}{2}) + pL_{2,j;x,y}
\mbox{}_{\sigma} \int_x^{\infty }F(\frac{z+y}{2})N(x,z)dz$$
$$ =  L_{2,j;x,y}F(\frac{x+y}{2}) + p\sum_l (a_l\sigma _x^l + (-1)^{jl} b_l \mbox{ }^1\sigma _z^l)
\mbox{}_{\sigma} \int_x^{\infty }F(\frac{z+y}{2})N(x,z)dz$$
$$= L_{2,j;x,y}F(\frac{x+y}{2}) + p \sum_l \{ [(a_l \mbox{ }^2\sigma _x^l
+ (-1)^{(j+1)l} b_l \mbox{ }^2\sigma _z^l)
\mbox{}_{\sigma} \int_x^{\infty }F(\frac{z+y}{2})N(x,z)dz]$$ $$  + a_l A_l(F;N)(x,y)+  (-1)^{jl} b_l B_l(F;N)(x,y) \} .$$
Imposing the condition
\par $(6)$ $L_{2;x,y} F(x,y)=0$, where \par $(7)$ $ L_{2;x,y} = L_{2,1;x,y} + L_{2,2;x,y} =
\sum_l ( 2 a_l \sigma _x^l + (1+ (-1)^l) b_l \sigma _y^l)$, we get the nonlinear
PDE with the help of Lemma 5 and Proposition 6
$$(8)\quad L_{2;x,y} K(x,y) - p \sum_l \{ 2 a_l \hat{A}_l(K;N)(x,y) +  (1+(-1)^l) b_l \hat{B}_l(K;N)(x,y) \} = 0 ,$$
since the center of the Cayley-Dickson algebra ${\cal A}_m$ is the real field $\bf R$ and so
the commutator of $bI$ and $\mbox{ }_{\sigma }\int $ is zero, $[bI, \mbox{ }_{\sigma }\int ] =0$,
also $(bI)(FK)= F(bIK)=bFK$, when $b$ is a real constant, where $I$ is the unit operator.
A solution of PDE $(8)$ can be found from the linear problem $(1,3,6,7)$ using
the integral operator $(4)$, where expressions for $\hat{A}_l$ and $\hat{B}_l$ are prescribed by Formulas
5$(3-7)$.
\par Making the variable change $y\mapsto - y$ one gets the PDO $\sigma _x + \sigma _y$
instead of $\sigma _x - \sigma _y$ and the corresponding changes in the PDO $L_2$.
Then Theorem 3 implies the following.

\par {\bf 8.1. Theorem.} {\it Let conditions of Theorem 3 and Example 8 be fulfilled,
then a solution of PDE $(8)$ is described by Formulas $(3,4,6)$, where PDOs $L_1$ and $L_2$ are provided by expressions $(1,7)$.}

\par {\bf 9. Example.} Consider now the generalization of PDOs from \S 5 with $k\ge 2$:
\par $(1)$ $L_1 =L_{1;x,y}= \sigma _x^k - \sigma _y^k$,
\par $(2)$ $L_{2,j} =L_{2,j;x,y}= \sum_l (a_l \sigma _x^{kl}  +   (-1)^{jkl} b_l \sigma _y^{kl})$,
\\ where $k$ is a natural number, $a_l, b_l \in {\bf R}$ for each $l$ when the Dirac type operator
$\sigma $ is over ${\cal A}_m$ with $m\ge 2$, the sum is finite or infinite, $j=1$ or $j=2$. Other suppositions are as in \S 8.
Let \par $(3)$ $L_{1;x,y}F(x,y) =0$ and
\par $(4)$ $K(x,y) = F(x,y) + p \mbox{}_{\sigma} \int_x^{\infty }F(z,y)N(x,z)dz$,
\\ where $N$ is expressed through $K$ by 1$(1,2)$.
\par Then we deduce the identities:
$$(5)\quad L_{2,j;x,y}K(x,y) = L_{2,j;x,y}F(x,y) + pL_{2,j;x,y}
\mbox{}_{\sigma} \int_x^{\infty }F(z,y)N(x,z)dz$$
$$ =  L_{2,j;x,y}F(x,y) + p\sum_l (a_l \sigma _x^{kl} + (-1)^{jkl} b_l \mbox{ }^1\sigma _z^{kl})
\mbox{}_{\sigma} \int_x^{\infty }F(z,y)N(x,z)dz$$
$$= L_{2,j;x,y}F(x,y) + p \sum_l \{ [(a_l \mbox{ }^2\sigma _x^{kl}
+ (-1)^{(j+1)kl} b_l \mbox{ }^2\sigma _z^l)
\mbox{}_{\sigma} \int_x^{\infty }F(z,y)N(x,z)dz]$$ $$  + a_l \hat{A}_{kl}(F;N)(x,y) + (-1)^{jkl} b_l \hat{B}_{kl}(F;N)(x,y) \} .$$
From the condition
\par $(6)$ $L_{2;x,y} F(x,y)=0$ with the PDO
\par $(7)$ $ L_{2;x,y} = L_{2,1;x,y} + L_{2,2;x,y} =
\sum_l ( 2 a_l \sigma _x^{kl} +  (1+ (-1)^{kl}) b_l \sigma _y^{kl})$ \\ we infer that
a function $K$ is a solution of the nonlinear PDE of the form:
$$(8)\quad L_{2;x,y} K(x,y) - p \sum_l \{ 2 a_l \hat{A}_{kl}(K;N)(x,y) +  (1+(-1)^{kl}) b_l \hat{B}_{kl}(K;N)(x,y) \} =0.$$
PDE $(8)$ can be resolved with the help of the linear problem $(1,3,6,7)$ and the integral operator $(4)$,
where terms $\hat{A}_l$ and $\hat{B}_l$ are given by Formulas 5$(3-7)$, a function $N$
satisfies conditions 1$(1-4)$ and either 2$(13)$ or 2$(14)$.
Thus due to Theorem 3 we have proved the following.

\par {\bf 9.1. Theorem.} {\it Let conditions of Theorem 3 and Example 9 be satisfied,
then a solution of PDE $(8)$ is provided by Formulas $(3,4,6)$, where PDOs $L_1$ and $L_2$ are given by $(1,7)$.}

\par {\bf 10. Example.} Let now the pair of PDOs be
\par $(1)$ $L_1 =L_{1;x,y}= \sigma _x^k + \sigma _y^k$,
\par $(2)$ $L_{2,j} =L_{2,j;x,y}= \sum_l (a_l \sigma _x^{kl}  + (-1)^{lj(k+1)} b_l \sigma _y^{kl})$,
\\ where $k$ is a natural number, $k\ge 2$, $K$, $F$, $U$ and $\sigma $ have the same meaning
as in \S \S 1 and 2, $a_l, b_l \in {\bf R}$ for each $l$ when $\sigma $ are over ${\cal A}_m$
with $m\ge 2$, $F\in Mat_{n\times n} ({\bf R})$ and $K\in Mat_{n\times n} ({\cal A}_r)$,
the sum is finite or infinite, $j=1$ or $j=2$.  Imposing the conditions
\par $(3)$ $L_{1;x,y}F(x,y) =0$ and
\par $(4)$ $L_{2;x,y} F(x,y)=0$ with \par $(5)$ $ L_{2;x,y} = L_{2,1;x,y} + L_{2,2;x,y} =
\sum_l ( 2 a_l \sigma _x^{kl} + (1+ (-1)^{l(k+1)}) b_l \sigma _y^{kl})$ and considering the integral transform
\par $(6)$ $K(x,y) = F(x,y) + p \mbox{}_{\sigma} \int_x^{\infty }F(z,y)N(x,z)dz,$
\\ where $N$ is related with $K$ by expressions 1$(1,2)$,
\\ we infer that
$$(7)\quad L_{2,j;x,y}K(x,y) = L_{2,j;x,y}F(x,y) + pL_{2,j;x,y}
\mbox{}_{\sigma} \int_x^{\infty }F(z,y)N(x,z)dz$$
$$ =  L_{2,j;x,y}F(x,y) + p\sum_l (a_l\sigma _x^{kl} + (-1)^{l(j+1)+jkl} b_l \mbox{ }^1\sigma _z^{kl})
\mbox{}_{\sigma} \int_x^{\infty }F(z,y)N(x,z)dz$$
$$= L_{2,j;x,y}F(x,y) + p \sum_l \{ [(a_l \mbox{ }^2\sigma _x^{kl}
+ (-1)^{l(j+1)(k+1)} b_l \mbox{ }^2\sigma _z^l)
\mbox{}_{\sigma} \int_x^{\infty }F(z,y)N(x,z)dz]$$ $$  + a_l A_{kl}(F;N)(x,y) + (-1)^{l(j+1)+jkl} b_l B_{kl}(F;N)(x,y) \} .$$
Thus in virtue of Lemma 5 and Proposition 6 a function $K$ satisfies the nonlinear PDE:
$$(8)\quad L_{2;x,y} K(x,y) - p \sum_l \{ 2 a_l \hat{A}_{kl}(K;N)(x,y) + ((-1)^l+(-1)^{kl}) b_l \hat{B}_{kl}(K;N)(x,y) \} =0.$$
The solution of the latter PDE reduces to the linear problem $(1,3-5)$ and using the integral operator $(6)$,
where terms $\hat{A}_l$ and $\hat{B}_l$ are given by Formulas 5$(3-7)$, a function $N$
is of the form 1$(1-4)$ and either 2$(13)$
or 2$(14)$ is fulfilled also. It is also possible to change the notation
$b_l\mapsto (-1)^lb_l$ in this example or $b_l\mapsto - b_l$ in \S 5.
Examples 7-10 correspond to different types of PDEs such as elliptic, hyperbolic and mixed types.
In view of Theorem 3 this implies the following.

\par {\bf 10.1. Theorem.} {\it If conditions of Theorem 3 and Example 10 are satisfied,
then a solution of PDE $(8)$ is given by Formulas $(3,4,6)$, where PDOs $L_1$ and $L_2$ are as in $(1,5)$.}

\par {\bf 10.2. Remark.} Transformation groups related with the quaternion skew field
are described in \cite{port69}. Automorphisms and derivations of the quaternion skew field
and the octonion algebra are contained in \cite{schaferb}, that of Lie algebras and groups
in \cite{gogros}.

\par {\bf 11. Example.} Consider now the term $N$ in the integral operator
\par $(1)$ $f(y)(g(x)K(x,y)) = F(x,y) + p \mbox{}_{\sigma} \int_x^{\infty }F(z,y)[f(z)(g(x)EK(x,z))]dz$ \\
with multiplier functions $f(z)$ and $g(x)$ satisfying definite conditions (see below),
where $F$, $K$ and $N(x,z) = f(z)(g(x)EK(x,z))$, $p$
have the meaning of the preceding paragraphs, $E$ is an operator fulfilling Conditions 1$(2,3)$ and either 2$(13)$
or 2$(14)$ also. Suppose that
\par $(2)$ $\sigma _z f(z) = \sum_j i_j \psi _j \partial f(z)/\partial z_{\xi (j)}= \lambda f(z)$ and
\par $(3)$ $\sigma _x g(x) = \sum_j i_j \psi _j \partial g(x)/\partial x_{\xi (j)}= \mu g(x)$, where
\par $\lambda = \sum_j i_j \psi _j \lambda _j$ and
\par $\mu = \sum_j i_j \psi _j \mu _j$ with $\lambda _j, \mu _j \in \bf R$ for each $j$.
We choose the functions $f(z) = C_1 \exp (\sum_j z_j \lambda _j)$ and $g(x) = C_2 \exp (\sum_j x_j \mu _j)$
satisfying PDEs $(2)$ and $(3)$ correspondingly, where
$C_1$ and $C_2$ are real non-zero constants, $x_j , z_j \in \bf R$, $x=\sum_j i_jx_j $, $~ x, z\in U$.
The first PDO we take as
\par $(4)$ $L_1= L_{1,x,y} = \sigma _x^k + s \sigma _y^k$, \\ where $k\ge 1$,
either $s=1$ or $s=-1$.
Then the condition
\par $(5)$ $L_{1,x,y} F(x,y)= 0$ is equivalent to
\par $(6)$ $\sigma _x^kF(x,y) = - s \sigma _y^kF(x,y).$ \par
Therefore we get from Proposition 2.5 with $N(x,z) = f(z)(g(x)EK(x,z))$ that
$$(7)\quad \sigma _y^{kl} \mbox{}_{\sigma} \int_x^{\infty }F(z,y)[f(z)(g(x)EK(x,z))]dz
= (- s)^l \mbox{ }^1\sigma _z^{kl} \mbox{}_{\sigma} \int_x^{\infty }F(z,y)[f(z)(g(x)EK(x,z))]dz$$
$$ = s^l (-1)^{l(k+1)} [\mbox{ }^2\sigma _z + \mbox{ }^4\sigma _z]^{kl} \mbox{}_{\sigma} \int_x^{\infty }F(z,y)[f(z)(g(x)EK(x,z))]dz$$
$$+ (-s)^l B_{kl}(F(z,y);[f(z)(g(x)EK(x,z))])|_{z=x},$$
where $F$ stands on the first place, $f$ on the second, $g$ on the third and $(EK)$ on the fourth place.
Then from $(2)$ and $(7)$ it follows that
$$(8)\quad \sigma _y^{kl} \mbox{}_{\sigma} \int_x^{\infty }F(z,y)[f(z)(g(x)EK(x,z))]dz= $$
$$s^l (-1)^{l(k+1)} [\mbox{ }^4\sigma _z + \lambda ]^{kl} \mbox{}_{\sigma} \int_x^{\infty }F(z,y)[f(z)(g(x)EK(x,z))]dz$$
$$+ (-s)^l B_{kl}(F(z,y);[f(z)(g(x)EK(x,z))])|_{z=x}.$$
Evaluation of the other integral with the help of Proposition 2.5 and Formula $(3)$ leads to:
$$(9)\quad \sigma _x^{kl} \mbox{}_{\sigma} \int_x^{\infty }F(z,y)[f(z)(g(x)EK(x,z))]dz
=$$  $$ [\mbox{ }^3\sigma _x + \mbox{ }^4\sigma _x]^{kl} \mbox{}_{\sigma} \int_x^{\infty }F(z,y)[f(z)(g(x)EK(x,z))]dz$$
$$+ A_{kl}(F(z,y);[f(z)(g(x)EK(x,z))])|_{z=x},$$
$$ = [\mbox{ }^4\sigma _x + \mu ]^{kl} \mbox{}_{\sigma} \int_x^{\infty }F(z,y)[f(z)(g(x)EK(x,z))]dz$$
$$+ A_{kl}(F(z,y);[f(z)(g(x)EK(x,z))])|_{z=x}.$$
Thus in this particular case PDEs of Examples 7-10 change. For example, PDE 10$(8)$ takes the form:
 $$(10)\quad \sum_l (2a_l (\sigma _x + \mu )^{kl}  + (1+(-1)^{l(k+1)}) b_l (\sigma _y + \lambda )^{kl})
K(x,y)$$  $$ - \frac{p}{f(y)g(x)}  \sum_l \{ 2 a_l \hat{A}_{kl}([f(y)(g(z) K(z,y))];[f(z)(g(x) E K(x,z))])|_{z=x} $$  $$ +  (1+(-1)^{l(k+1)}) b_l \hat{B}_{kl}([f(y)(g(z) K(z,y))];[f(z)(g(x) E K(x,z))])|_{z=x} \} =0  ,$$
when $F \in Mat_{n\times n}({\bf R})$ and $K\in Mat_{n\times n}({\cal A}_r)$ with $2\le r$ and
$E$ is the right linear operator over ${\cal A}_r$, since the operator $E$ satisfies Conditions
1$(2,3)$ and either 2$(13)$
or 2$(14)$;
the functions $f(y)$ and $g(x)$ have values in ${\bf R}\setminus \{ 0 \} $ for each
$x, y \in U$, whilst ${\bf R}$ is the center of the Cayley-Dickson algebra.
Analogous changes will be in Examples 7-9.

\par {\bf 12. Example.}
Let the non-commutative integral operator be
\par $(1)$ $K(x,y) = F(x,y) + {\sf B}_xK(x,y)$ with
$$(2)\quad {\sf B}_xK(x,y)=p \mbox{}_{\sigma} \int_x^{\infty }F(z,y)N(x,z,y)dz,$$
where $F$, $K$, $N(x,z,y)$ are as in Proposition 2.5 and Theorem 3.3,
$F\in Mat_{n\times n}({\bf R})$, $K$ and $N$ are in $Mat_{n\times n}({\cal A}_m)$,
while $p$ is a sufficiently small non-zero real parameter, $N$ is an operator function right linear in $K$
as in \S 1. Put
\par $(3)$ $N(x,z,y) = E_yK(x,z)$ for every $x$, $y$ and $z$ in $U$, \\ where $[E_y,L_j]=0$ for each $j=1,...,k_0$ and $y\in U$, $E=E_y$ may depend
on the variable $y\in U$ also, $E$ is an operator satisfying Conditions 1$(2,3)$ and either 2$(13)$
or 2$(14)$, $m\ge 2$.
Choose two PDO
\par $(4)$ $L_1= L_{1,x,y} = \sigma _x + \sigma _y$,
\par $(5)$ $L_2 = L_{2,x,y} = (\sum_l a_l \sigma _x^l) + s \sigma _y$,
\\ where $s\in {\bf R}$, $s$ is a non-zero real constant,
$a_l \in {\bf R}$ for each $l$ when the Dirac type operator
$\sigma $ is over ${\cal A}_m$ with $m\ge 2$. We impose the conditions:
\par $(6)$ $ L_{j,x,y} F(x,y) = 0$ for $j=1$ and $j=2$, for all $x, y \in U$.
Then it is possible to write $F(z,y) = F(\frac{z-y}{2})$.
Applying the PDO $L_2$ to both sides of $(1)$ and using $(2)$, Proposition 2.5
and Conditions $(3-6)$ we deduce that
$$(7)\quad L_{2,x,y} K(x,y) =  p \{ \sum_l a_l \sigma _x^l - s\mbox{ }^1\sigma _z + s\mbox{ }^2\sigma _y \} \mbox{}_{\sigma} \int_x^{\infty }F(z,y)N(x,z,y)dz$$
$$ =  [p \{ \sum_l a_l \mbox{ }^2\sigma _x^l + s \mbox{ }^2\sigma _z + s\mbox{ }^2\sigma _y \}  \mbox{}_{\sigma} \int_x^{\infty }F(z,y)N(x,z,y)dz]$$
$$+ [p \sum_l a_l A_l (F;N)(x,y)] - ps B_1(F;N)(x,y).$$ For sufficiently small non-zero real values of $p$ the operator
$I-{\sf B}_x$ is invertible and hence Equality $(7)$, Lemma 5 and Proposition 6
imply that $K$ satisfies the nonlinear partial integro-differential equation:
$$(8)\quad L_{2,x,y} K(x,y) -  [p \sum_l a_l \hat{A}_l (K;N)(x,y)] - ps K(x,y)N(x,x,y)$$ $$ - ps \mbox{}_{\sigma} \int_x^{\infty }K(z,y)\sigma _y N(x,z,y)dz=0 .$$
Using Theorem 3 we deduce the following.

\par {\bf 12.1. Theorem.} {\it  A solution of PIDE $(8)$ is described by $(1,2,3,6)$, where PDOs $L_1$ and $L_2$ are given by $(4,5)$,
provided that conditions of Theorem 3 and Example 12 are satisfied.}

\par {\bf 13. Example.} Suppose that functions $K$ and $F$ are related by Equations
12$(1,2)$ and take two PDOs
\par $(1)$ $L_{1,x,y} = \sigma _x - \sigma _y$ and
\par $(2)$ $L_{2,x,y} = \Delta _x + s \Delta _y$,
\\ where the coefficient $\psi _0$ is null in $\sigma $ and hence the Laplace operator is
expressed as $\Delta = - \sigma ^2$, while $s\in {\bf R}\setminus \{ 0 \} $.
Now we take a function $N$ in the  form
\par $(3)$ $N(x,z,y) = E K(x,ay+bz)$, \\ where $a$ and $b$ real parameters
to be calculated below such that $a^2+b^2>0$, $ ~ b$ is non-zero. Then from the conditions
\par $(4)$ $L_{j,x,y}F(x,y)=0$ for $j=1$ and $j=2$, \\
Proposition 2.5 and Corollary 2.6 it follows that
$$(5)\quad L_{2,x,y}K(x,y) = - p (\sigma _x^2 + s \sigma _y^2)
\mbox{}_{\sigma} \int_x^{\infty }F(z,y)E K(x,ay+bz) dz$$
$$ = p (\mbox{ }^2\Delta _x - s (\mbox{ }^1\sigma _y + \mbox{ }^2\sigma _y)^2) \mbox{}_{\sigma} \int_x^{\infty }F(z,y)E K(x,ay+bz) dz - p A_2(F(z,y), E K(x,ay+bz) )|_{z=x}$$ and
$$(6)\quad (\mbox{ }^1\sigma _y + \mbox{ }^2\sigma _y)^2 \mbox{}_{\sigma} \int_x^{\infty }F(z,y)E K(x,ay+bz) dz
= $$ $$[\mbox{ }^1\sigma _z^2+  ab^{-1}\mbox{ }^1\sigma _z \mbox{ }^2\sigma _z + ab^{-1} \mbox{ }^2\sigma _z \mbox{ }^1\sigma _z + a^2b^{-2} \mbox{ }^2\sigma _z^2] \mbox{}_{\sigma} \int_x^{\infty }F(z,y)E K(x,ay+bz) dz$$
$$ = [ab^{-1}(\mbox{ }^1\sigma _z + \mbox{ }^2\sigma _z)^2+  (1-ab^{-1})\mbox{ }^1\sigma _z^2 +
(a^2b^{-2} - ab^{-1}) \mbox{ }^2\sigma _z^2] \mbox{}_{\sigma} \int_x^{\infty }F(z,y)E K(x,ay+bz) dz$$
$$ = (1-ab^{-1})^2\mbox{ }^2\sigma _z^2 \mbox{}_{\sigma} \int_x^{\infty }F(z,y)E K(x,ay+bz) dz $$
$$ - ab^{-1}[\sigma _z (F(z,y) E K(x,ay+bz) )]|_{z=x} + (1-ab^{-1})B_2(F,EK)(x,y)$$
$$ = p^{-1} (1-ab^{-1})^2b^2{\sf B}_x(\sigma _y^2K(x,y))$$ $$ - ab^{-1}[\sigma _z (F(z,y) E K(x,ay+bz) )]|_{z=x} + (1-ab^{-1})B_2(F,EK)(x,y),$$
since $$\sigma _z^2 \mbox{}_{\sigma} \int_x^{\infty }F(z,y)N(x,z,y) dz =
 \mbox{}_{\sigma} \int_x^{\infty }\sigma _z^2[F(z,y)N(x,z,y)] dz$$ $$ =
 - \sigma _z[(F(z,y) E K(x,ay+bz) )]|_{z=x} .$$
Then Identities $(5,6)$ imply that
$$(7)\quad L_{2,x,y}K(x,y) = {\sf B}_x [L_{2,x,y}K(x,y)] - p A_2(F(z,y), E K(x,ay+bz) )|_{z=x}$$
$$ +ps ab^{-1}[\sigma _z (F(z,y) E K(x,ay+bz) )]|_{z=x} - ps (1-ab^{-1})B_2(F,EK)(x,y)$$
when $(b-a)^2=1$ and $b$ is non-zero, that is either $a=b+1$ or $a=b-1$.
In virtue of Lemma 5 and Proposition 6 this gives the nonlinear PDE for $K$:
$$(8)\quad L_{2,x,y}K(x,y) + p \hat{A}_2(K(z,y), E K(x,ay+bz) )|_{z=x}$$
$$- ab^{-1}ps[\sigma _z (K(z,y) E K(x,ay+bz) )]|_{z=x} + (1-ab^{-1})ps\hat{B}_2(K(z,y),EK(x,ay+bz))|_{z=x} = 0 .$$
From Theorem 3 we infer the following.

\par {\bf 13.1. Theorem.} {\it  A solution of PDE $(8)$ is given by $(3,4)$ and 12$(1,2)$, where PDOs $L_1$ and $L_2$ are
prescribed by $(1,2)$, whenever conditions of Theorem 3 and Example 13 are satisfied.}

\par {\bf 14. Nonlinear PDE with parabolic terms.}
Let \par $(1)$ $\partial _t := \sum_{k=1}^v  \partial /\partial t_k$ \\ be the first order
PDO, where $t_1,...,t_v$ are real variables independent of other variables
$x, y, z\in U$, $ ~ t=(t_1,...,t_v)\in W$, $W := \{ t\in {\bf R}^v: ~\forall k=1,...,v ~ 0\le t_k< T_k \} $, where
$T_k$ is a constant, $0<T_k\le \infty $ for each $k$.
\par Suppose that \par $(2)$ $F$ and $K$ are continuously differentiable functions by $t_k$ for each $k$ so that
$\mbox{ }_{\sigma }\int _x^{\infty } F(z,y)N(x,z,y)dz $ converges for some $t\in W$
and
\par $(3)$ the integrals $\mbox{ }_{\sigma }\int _x^{\infty } (\partial_t F(z,y))N(x,z,y)dz$ and
$\mbox{ }_{\sigma }\int_x^{\infty } F(z,y)(\partial _tN(x,z,y))dz $ converge uniformly on $W$
in the parameter $t$.
\par In virtue of the theorem about differentiation of an improper integral by a parameter
the equality is valid: $$(4)\quad \partial _t \mbox{ }_{\sigma }\int _x^{\infty } F(z,y)N(x,z,y)dz =$$
$$\mbox{ }_{\sigma }\int _x^{\infty } (\partial_t F(z,y))N(x,z,y)dz + \mbox{ }_{\sigma }\int_x^{\infty } F(z,y)(\partial _tN(x,z,y))dz .$$
Using $(4)$ the commutator
$(I- {\sf A}_xE_y)((\partial_t + L_s)f)- (\partial_t + L_s)[(I- {\sf A}_xE_y)f]$ can be calculated, when
there is possible to evaluate the commutator $(I- {\sf A}_xE_y)(L_sf)-L_s[(I- {\sf A}_xE_y)f]=R_s(f)$ for suitable functions $f$ and a PDO $L_s=L_{s,x,y}$ (see also \S 2).
\par For solution of nonlinear PDEs or PIDEs also the following will be useful
for integral operators of the form $\mbox{ }_{\sigma }\int_x^{\infty } N(x,z,y) K(x,z)dz$.

\par {\bf 14.1. Example.} Let a PDO be
\par $(1)$ $L_1=\partial _t + \sum_l a_l(\sigma _x^l +
(-1)^{l+1}\sigma _y^l)$ and let \par $(2)$ $N(x,z,y)=E_yK(x,z)$,
\\ where $a_l\in \bf R$ for all $l=0,1,2,...$, so that conditions
1$(2,3)$ and 2$(1,5)$ are fulfilled, $F\in Mat_{n\times n}({\bf
R})$, $K\in Mat_{n\times n}({\cal A}_r)$, $2\le r$, the first order
PDO $\sigma $ is over the Cayley-Dickson algebra ${\cal A}_r$ (see
\S \S 1, 2 and 14). Then
\par $(\partial _t+ \sum_l a_l (-1)^{l+1}\sigma _y^l)F(x,y) = -
(\sum_l a_l\sigma _x^l)F(x,y)$\\ for all $x, y \in U$. Therefore we
infer from Proposition 2.5 that $$(3)\quad L_1K(x,y) = p (\mbox{
}^2\partial _t + \sum_l a_l(\sigma _x^l - \mbox{ }^1\sigma _z^l))
\mbox{ }_{\sigma }\int_x^{\infty } F(z,y)E_yN(x,z)dz$$
$$=p (\mbox{ }^2\partial _t + \sum_l
a_l(\mbox{ }^2\sigma _x^l + (-1)^{l+1}\mbox{ }^2\sigma _z^l)) \mbox{
}_{\sigma }\int_x^{\infty } F(z,y)E_yN(x,z)dz$$
$$+p\sum_l a_l(A_l(F;N)(x,y)- B_l(F;N)(x,y)).$$ Hence we deduce a
nonlinear PDE
$$(4)\quad L_1K(x,y) - p\sum_l a_l(\hat{A}_l(K;EK)(x,y)- \hat{B}_l(K;EK)
(x,y))=0$$ according to Lemma 5 and Proposition 6. Its solution
reduces to the linear problem 2$(2,5)$. We mention that PDE $(4)$
corresponds to some kinds of Sobolev type nonlinear PDEs.

\par {\bf 14.2. Remark.} Suppose that $L$ and $S$ are PDOs and
functions $f: (a,b)\times U^m\to {\cal A}_r$ and $g: (a,b)\times
U^m\to {\cal A}_r$ are in the domains of operators $\exp (tL)$ and
$S$ correspondingly, where $t$ is a real parameter, $t\in (a,b)$,
$a<b$, $2\le r$, where PDOs $L$ and $S$ are by variables in $U^m$,
$m\in {\bf N}$. If they satisfy the PDE
\par $(1)$ $\exp (tL)f(t,x_1,...,x_m)=Sg(t,x_1,...,x_m)$\\
for all $t\in (a,b)$ and $x_1,...,x_m\in U$, then $$ \frac{\partial
\exp (tL)f(t,x_1,...,x_m)}{\partial t} = \exp (tL)(\frac{\partial
}{\partial t} + L)f(t,x_1,...,x_m)$$ $$= \exp (tL)(\frac{\partial
}{\partial t} + L)\exp (-tL)Sg(t,x_1,...,x_m)= \frac{\partial
Sg(t,x_1,...,x_m)}{\partial t} ,$$ consequently,
$$(2)\quad(\frac{\partial }{\partial t} + L)f(t,x_1,...,x_m)= \exp
(-tL)\frac{\partial Sg(t,x_1,...,x_m)}{\partial t}.$$ The latter
also may be helpful for solutions of nonlinear PDEs with parabolic
terms.

\par {\bf 14.3. Generalized approach.} Let $L_1,...,L_k$ and
$S_1,...,S_k$ be PDOs which are polynomials or series of $\sigma _x$
and $\sigma _y$ so that \par $(1)$ $[L_j,S_j]=0$ \\ for each
$j=1,...,k$, where $x$ and $y$ are in a domain $U$ in the
Cayley-Dickson algebra ${\cal A}_r$, $2\le r$ (see \S 2.3). Instead
of the conditions $L_jF=0$ it is
possible to consider more generally \par $(2)$ $L_jF=G_j$, \\
where $G_j$ are some functions known or defined by some relations,
while functions $F$, $G_j$ and $K$ may also depend on a parameter
$t\in W$ (see \S 14) so that $F\in Mat_{n\times n}({\bf R})$, $G_j$
for all $j$ and $K$ have values in $Mat_{n\times n}({\cal A}_r)$. It
is also supposed that $F$ and $K$ are related by the integral
equation 2$(1)$ and $N(x,y,z)=E_yK(x,z)$ and Conditions 1$(2,3)$ are
satisfied. In particular, if \par $(3)$ $G_j=L_j(I+S_j)K$, then a
solution of the linear system of PIDEs
\par $(4)$ $L_jF(x,y)=L_j(I+S_j)K(x,y)$ and
\par $(5)$ $(I-{\sf A}_xE)K(x,y)=F(x,y)$ \\
would also be a solution of nonlinear PIDEs
\par $(6)$ $S_jL_jK(x,y)+ M_j(K)=0$, \\
where $M_j$ corresponds to $L_j$ for each $j=1,...,k_0$ with $1\le
k_0\le k$ as in \S 2. Thus this generalizes PIDEs 2$(12)$.
Particularly, taking $S_j=\partial _t$ we get that Condition $(1)$
is valid. Therefore, the technique presented in \S 7-14.3
encompasses some kinds of nonlinear Sobolev type PDEs as well.

\par {\bf 15. Proposition.} {\it Let $$(1)\quad \lim_{z\to \infty } ~ \mbox{}^1\sigma _z^k ~ \mbox{}^2\sigma
_x^s ~ \mbox{}^2\sigma _z^nN(x,z,y)K(x,z)=0$$ for each $x,
y$ in a domain $U$ satisfying Conditions 2.1$(D1,D2)$ with $\infty
\in U$ and every non-negative integers $0\le k, s, n\in {\bf Z}$
such that $k+s+n\le m $. Suppose also that $\mbox{}_{\sigma }
\int_x^{\infty } \partial ^{\alpha }_x\partial ^{\beta }_y \partial
^{\omega }_z [N(x,z,y)K(x,z)] dz $ converges uniformly by
parameters $x, y $ on each compact subset $W\subset U\subset {\cal
A}_r^2$ for each $|\alpha |+|\beta |+|\omega |\le m$, where $\alpha
=(\alpha _0,...,\alpha _{2^r-1})$, $|\alpha |=\alpha _0+...+\alpha
_{2^r-1}$, $\partial ^{\alpha }_x=\partial ^{|\alpha |}/\partial
x_0^{\alpha _0} ...\partial x_{2^r-1}^{\alpha _{2^r-1}}$, where
$N\in C^m(U^3,Mat_{n\times n}({\cal A}_r))$ and $K\in C^m(U^2,Mat_{n\times n}({\cal A}_r))$.
Then
$$(2)\quad  \sigma _x^m\mbox{ }_{\sigma }\int_x^{\infty } N(x,z,y) K(x,z) dz =$$
$$(\mbox{ }^1\sigma _x + \mbox{ }^2\sigma _x)^{m} \mbox{ }_{\sigma }\int_x^{\infty } N(x,z,y) K(x,z) dz
+ \tilde{A}_m(N;K)(x,y),$$
where
$$(3)\quad \tilde{A}_m(N;K)(x,y) = - \sum_{j=0}^{m-1} \sigma _x^j \{ [ \sigma _x^{m-j-1}N(x,z,y)K(x,z)]|_{z=x} \} $$
for each $m\ge 1$, $ \sigma _x^0=I$. Moreover,
$$(4)\quad  \mbox{ }^1\sigma _z^m \mbox{ }_{\sigma }\int_x^{\infty } N(x,z,y) K(x,z) dz =$$
$$ (- \mbox{ }^2\sigma _z)^m \mbox{ }_{\sigma }\int_x^{\infty } N(x,z,y) K(x,z) dz
+\tilde{B}_m(N;K)(x,y), $$ where
$\tilde{B}_m(N;K)(x,y) = \tilde{B}_m(N(x,z,y);K(x,z))|_{z=x}$,
$$(5)\quad \tilde{B}_m(N(x,z,y);K(x,z)) = [\sum_{k=0}^{m-1} (-1)^{k+1} \mbox{ }^1\sigma _z^{m-k-1}\mbox{ }^2\sigma _z^{k}]
N(x,z,y) K(x,z).$$  }
\par {\bf Proof.} For $m=1$ we infer that
$$(6)\quad \sigma _x\mbox{ }_{\sigma }\int_x^{\infty } N(x,z,y) K(x,z) dz =$$
$$(\mbox{ }^1\sigma _x + \mbox{ }^2\sigma _x)  \mbox{ }_{\sigma }\int_x^{\infty } N(x,z,y) K(x,z) dz - N(x,x,y)K(x,x)$$
and put $\tilde{A}_1(N;K)(x,y) = \tilde{A}_1(N(x,z,y);K(x,z))|_{z=x}$, where
\par $(7)$ $\tilde{A}_1(N(x,z,y);K(x,z)) = - N(x,z,y) K(x,z)$. \par Then we deduce by induction that
$$(8)\quad \sigma _x^m\mbox{ }_{\sigma }\int_x^{\infty } N(x,z,y) K(x,z) dz =$$
$$\sigma _x[(\mbox{ }^1\sigma _x + \mbox{ }^2\sigma _x)^{m-1} \mbox{ }_{\sigma }\int_x^{\infty } N(x,z,y) K(x,z) dz
+ \tilde{A}_{m-1}(N(x,z,y);K(x,z))|_{z=x}]=$$
$$(\mbox{ }^1\sigma _x + \mbox{ }^2\sigma _x)^{m} \mbox{ }_{\sigma }\int_x^{\infty } N(x,z,y) K(x,z) dz
+ \tilde{A}_m(N;K)(x,y)$$
with the convention that the operator $\mbox{ }_{\sigma }\int_x^{\infty }$ stands in the zero position,
$N$ in the first and $K$ in the second positions correspondingly,
where
$$(9)\quad \tilde{A}_m(N;K)(x,y) = \sigma _x \tilde{A}_{m-1}(N;K)(x,y) - [ \sigma _x^{m-1}N(x,z,y)K(x,z)]|_{z=x}$$
for each $m\ge 2$. Therefore, by induction we deduce that
\par $(10)$ $ \tilde{A}_m(N;K)(x,y) = - [ \sigma _x^{m-1}N(x,z,y)K(x,z)]|_{z=x} -$ \par $ \sigma _x \{ [ \sigma _x^{m-2}N(x,z,y)K(x,z)]|_{z=x}
\} -...- $\par $\sigma _x^{m-2} \{ [ \sigma _xN(x,z,y)K(x,z)]|_{z=x} \} - \sigma _x^{m-1} N(x,x,y)K(x,x)$.
\par Then we infer:
$$(11)\quad \mbox{ }^1\sigma _z \mbox{ }_{\sigma }\int_x^{\infty } N(x,z,y) K(x,z) dz =$$
$$ - \mbox{ }^2\sigma _z \mbox{ }_{\sigma }\int_x^{\infty } N(x,z,y) K(x,z) dz
+\tilde{B}_1(N;K)(x,y), $$ where
$\tilde{B}_1(N;K)(x,y) = \tilde{B}_1(N(x,z,y);K(x,z))|_{z=x}$,
\par $\tilde{B}_1(N(x,z,y);K(x,z)) = - N(x,z,y) K(x,z).$
\par Therefore applying the operator $\mbox{ }^1\sigma _z$ by induction we get the formulas
$$(12)\quad \mbox{ }^1\sigma _z^m \mbox{ }_{\sigma }\int_x^{\infty } N(x,z,y) K(x,z) dz =$$
$$ \mbox{ }^1\sigma _z^{m-1}[ - \mbox{ }^2\sigma _z \mbox{ }_{\sigma }\int_x^{\infty } N(x,z,y) K(x,z) dz]
-   [\mbox{ }^1\sigma _z^{m-1}N(x,z,y) K(x,z)]|_{z=x}$$
$$=...= (- \mbox{ }^2\sigma _z)^m \mbox{ }_{\sigma }\int_x^{\infty } N(x,z,y) K(x,z) dz
+\tilde{B}_m(N;K)(x,y), $$ where
$\tilde{B}_m(N;K)(x,y) = \tilde{B}_m(N(x,z,y);K(x,z))|_{z=x}$,
$$(13) \quad \tilde{B}_m(N(x,z,y);K(x,z)) = [\sum_{k=0}^{m-1} (-1)^{k+1} \mbox{ }^1\sigma _z^{m-k-1}\mbox{ }^2\sigma _z^{k}]
N(x,z,y) K(x,z)$$ $$ = - (-\mbox{ }^2\sigma _z)^{m-1} N(x,z,y)K(x,z) + \mbox{ }^1\sigma _z \tilde{B}_{m-1}(N(x,z,y);K(x,z)).$$

\par {\bf 16. Theorem.} {\it Let $\{ L_s: ~ s=1,...,k_0 \} $ be a set of PDOs which are polynomials
$\Omega _s(\mbox{ }_1\sigma _x, \mbox{ }_2\sigma _y)$ over ${\cal A}_r$ or $\bf R$.
Let also $G$ be the family of all operators $E=BST_g$ satisfying the condition $[L_s,E]=0$ for each $s=1,...,k_0$, where $B\in SL_n({\bf R})$, $~S\in Aut (Mat_{n\times n}({\cal A}_r))$,
$~g\in Diff^{\infty }(U)$, $T_g$ is prescribed by Formula 1$(4)$, $~ \mbox{ }_1\sigma _x$ and $\mbox{ }_2\sigma _y$ are over the Cayley-Dickson algebra ${\cal A}_r$, $r\ge 2$. Then the family $G$ forms the group
and there exists an embedding of $G$ into $SL_n({\bf R})\times Aut (Mat_{n\times n}({\cal A}_r))\times Diff^{\infty }(U)$.}
\par {\bf Proof.} The composition (set theoretic) in the family $G$ of the aforementioned operators
is associative. Then the inverse $E^{-1} = T_g^{-1}S^{-1}B^{-1}$ of $E=BST_g$ exists, since $B$, $S$ and $T_g$
are invertible for every $B\in SL_n({\bf R})$, $~S\in Aut (Mat_{n\times n}({\cal A}_r))$ and
$~g\in Diff^{\infty }(U)$ so that $T_g^{-1} = T_{g^{-1}}$. On the other hand, the identity
$E^{-1} [L_s,E]E^{-1}  = - [L_s,E^{-1}]$ is valid. Thus the equality $[L_s,E]=0$ implies that $L_s$ and $E^{-1}$
commute, $[L_s,E^{-1}] = 0$, as well. Therefore, from $E\in G$ the inclusion $E^{-1}\in G$ follows.
The identity $[L_s,E_1E_2] = [L_s,E_1]E_2 + E_1 [L_s,E_2]$ implies that $E_1E_2 \in G$ whenever
$E_1\in G$ and $E_2\in G$. Thus the family $G$ has the group structure.
There exists the bijective correspondence between diffeomorphisms $g\in Diff^{\infty }(U)$
and operators $T_g$ acting on functions defined on $U$ with values in $Mat_{n\times n}({\cal A}_r)$ according to Formula 1$(4)$.
Each element $E$ in $G$ is of the form $E=BST_g$, where $B\in SL_n({\bf R})$, $~S\in Aut (Mat_{n\times n}({\cal A}_r))$,
$~g\in Diff^{\infty }(U)$, consequently, an embedding $ \omega : G\hookrightarrow SL_n({\bf R})\times Aut (Mat_{n\times n}({\cal A}_r))\times Diff^{\infty }(U)$ exists.

\section {Nonlinear PDEs used in hydrodynamics.}
\par {\bf 1. Remark.} In the previous article \cite{lucveleq2013} vector hydrodynamical
PDEs were investigated. Using results of Sections 2 and 3 we generalize the approach
using transformations of functions by operators $E$ of the form 3.1$(2)$. It permits to consider
other PDEs and study the symmetry of solutions.
\par {\bf 2. Example. Generalized Korteweg-de-Vries' type PDE.}
Let \par $(1)$ $N(x,z,y) = EK(x,z)$ as in 3.1$(1)$ and let ${\sf A}_x$
be given by 3.2$(3)$, where $E$ satisfies
conditions 3.1$(2,4)$. Foliations of a domain $U$, operators $L_s$ and the cases of
$F$ and $K$ are the same as in \cite{lucveleq2013}:
\par $(2)$ $L_1=\mbox{}_1\sigma _x^2- ~ \mbox{}_2\sigma _y^2$ and
\par $(3)$ $L_2=\mbox{}_3\sigma _t + ~ \mbox{}_1\sigma ^3_x+
~ 3 ~ \mbox{}_2\sigma _y ~ \mbox{}_1\sigma ^2_x + ~ 3 ~
\mbox{}_2\sigma ^2_y ~ \mbox{}_1\sigma _x + ~ \mbox{}_2\sigma ^3_y$, \\
where $\mbox{}_1\psi _0= ~\mbox{}_2\psi _0 =0$,
\par $(4)$ $L_1F=0$ and $L_{2,j}F=0$ for each $j=0,...,2^{r-1}$.
Then equations from example 4.2 \cite{lucveleq2013} take the following form.
In view of 2.6$(7)$ and Proposition 2.5 above PDE $(4.19)$ \cite{lucveleq2013} transforms into:
\par $(5)$ $(\mbox{}_1\sigma _x^2- ~ \mbox{}_2\sigma _y^2)
K(x,y)+2p K(x,y)[~\mbox{}_1\sigma _xEK(x,x)]=0$.
\par Putting \par $(6)$ $u(x)=2~\mbox{}_1\sigma _x EK(x,x)$ \\ over the
quaternion skew field ${\bf H}={\cal A}_2$ and
substituting $K(x,y)=\Phi (x,k)\exp (J Re (ky))$ into $(4)$, we
get Schr\"odinger's equation:
\par $(7)$ $~\mbox{}_1\sigma _x^2\Phi (x,k) + \Phi (x,k)(pu -
\sum_jk_j^2~\mbox{}_2\psi _j^2)=0$, \\
where $k\in {\bf H}$, the generator $J$ commutes with $i_0,...,i_{2^r-1}$.
There is supposed that functions $F$ and $K$ may depend
on $t$. Then in formulas $(4.24-28)$ \cite{lucveleq2013} $K$ changes into $EK$, while
$(4.30)$ \cite{lucveleq2013} due to $(5)$, 3.1$(2-5)$ and 3.2$(1,2)$ transforms into:
\par $(8)$ $ -2 {\sf p} ~\mbox{}_{\mbox{}_1\sigma } \int _x^{\infty }
F(z,y) [ K(x,y) (\mbox{}_1\sigma _x K(x,x))]dz =
[K(x,y)-F(x,y)]ST_gu(x)$, \\
since $E(Ku) = B((ST_g)(Ku))=B(ST_gK)(ST_gu)= (EK)(ST_gu)$ (see $u$ in $(6)$).
Therefore $(4.29)$ \cite{lucveleq2013} changes into
\par $(9)$ $I_2= -3 ~ \mbox{}_2\sigma _y [K(x,y) - F(x,y)]
(ST_gu(x))$\par $ + ~ 3p ~ \mbox{}_2\sigma _y [\mbox{}_1A_2(F,E
K)(x,y) - ~\mbox{}_1B_2(F,EK)(x,y)]$.
\par Hence from Formulas $(5,8,9)$ it follows that
\par $(10)\quad (\mbox{}_3\sigma _t + ~ \mbox{}_1\sigma ^3_x+
~ 3 ~ \mbox{}_2\sigma _y ~ \mbox{}_1\sigma ^2_x + ~ 3 ~
\mbox{}_2\sigma ^2_y ~ \mbox{}_1\sigma _x + ~ \mbox{}_2\sigma
^3_y) K(x,y) + 3\mbox{}_2^1\sigma _y[K(x,y)ST_gu(x)]=$
\par $p (\mbox{}_3^2\sigma _t + ~ \mbox{}_1^2\sigma ^3_x+ ~ 3 ~
\mbox{}_1^2\sigma _z ~ \mbox{}_1^2\sigma ^2_x + ~ 3 ~
\mbox{}_1^2\sigma ^2_z ~ \mbox{}_1^2\sigma _x + ~ \mbox{}_1^2\sigma
^3_z)~ \mbox{}_{\mbox{}_1\sigma } \int_x^{\infty } F(z,y)E
K(x,z)dz$\par $ +3p~ \mbox{}_{\mbox{}_1\sigma }
\int_x^{\infty } F(z,y)[\mbox{}^1_2\sigma _z E(
K(x,z)u(x))]dz  +T$, where \par $T= p ~\mbox{}_1A_3(
F,EK)(x,y) - p ~\mbox{}_1B_3(F,E K)(x,y)$\par $ + ~3
 p ~ \mbox{}_2 \sigma _y [\mbox{}_1A_2(F,EK)(x,y) -
~ \mbox{}_1B_2(F,EK)(x,y)] +$\par $ ~ 3 p ~
\mbox{}_2^1\sigma _y[F(x,y)ST_gu(x)] + 3 p ~(~\mbox{}_1^2\sigma _z
~ \mbox{}_1^2\sigma _x + ~ \mbox{}_1^2\sigma _z^2 - ~
\mbox{}_1^1\sigma _x ~ \mbox{}_1^2\sigma _x - \mbox{}_1^1\sigma _x ~
\mbox{}_1^2\sigma _z) [F(x,y)EK(x,z)]|_{z=x}$.
\par Then using $(5,8-10)$ we infer that
\par $(11)\quad T=- p (3 ~\mbox{}^2_1\sigma ^2_x + ~ \mbox{}^2_1\sigma _x ~ \mbox{}^2_1\sigma
_z+ 2~ \mbox{}^2_1\sigma _z ~ \mbox{}^2_1\sigma _x) [
F(x,y)EK(x,z)]|_{z=x}$\par $ - p (2 ~\mbox{}^1_1\sigma
_x ~\mbox{}^2_1\sigma _x + ~\mbox{}^2_1\sigma _x ~\mbox{}^1_1\sigma
_x) [F(x,y)E K(x,x)]$\par $+ 3(1- p) ~\mbox{}_2\sigma _y [F(x,y)ST_gu(x)] + 3p ~(~\mbox{}_1^2\sigma _z ~
\mbox{}_1^2\sigma _x + ~ \mbox{}_1^2\sigma _z^2 - ~
\mbox{}_1^1\sigma _x ~ \mbox{}_1^2\sigma _x - \mbox{}_1^1\sigma _x ~
\mbox{}_1^2\sigma _z) [F(x,y)EK(x,z)]|_{z=x} $
\par $= - 3 p ~\mbox{}_1^1\sigma _x[F(x,y)u(x)] -
3p (\mbox{}_1^2\sigma _x^2- ~ \mbox{}_1^2\sigma _z^2)[F(x,y)EK(x,z)]|_{z=x}$\par $ + 3(1-p) ~\mbox{}_2\sigma _y [
F(x,y)ST_gu(x)] $
\par $+ p [ \mbox{}^2_1\sigma _z, ~ \mbox{}^2_1\sigma _x]
[F(x,y) EK(x,z)]|_{z=x} + p [ \mbox{}^1_1\sigma _x,
~ \mbox{}^2_1\sigma _x] [F(x,y) EK(x,x)]$
\par $= - 3p ~\mbox{}_1^1\sigma _x[F(x,y)u(x)] + 3p F(x,y)[
E(K(x,x)u(x))]+ 3(1- p) ~\mbox{}_2\sigma _y [F(x,y)ST_gu(x)] $
\par $+ p [ \mbox{}^2_1\sigma _z, ~ \mbox{}^2_1\sigma _x]
[F(x,y) EK(x,z)]|_{z=x} + p [ \mbox{}^1_1\sigma _x,
~ \mbox{}^2_1\sigma _x] [F(x,y) EK(x,x)]$
\par $=- 3p ~\mbox{}_1^1\sigma _x[K(x,y)u(x)] + 3p^2 ~\mbox{}_1\sigma
_x\{ (\mbox{}_{\mbox{}_1\sigma }\int_x^{\infty } F(z,y)E
K(x,z)dz)ST_gu(\eta ) \}|_{\eta =x} + 3 p F(x,y)[E(
K(x,x)u(x))]+ 3(1- p) ~\mbox{}_2\sigma _y [
F(x,y)ST_gu(x)] $ \par $ + p [ \mbox{}^2_1\sigma _z, ~ \mbox{}^2_1\sigma
_x] [F(x,y) EK(x,z)]|_{z=x} + p [ \mbox{}^1_1\sigma
_x, ~ \mbox{}^2_1\sigma _x] [F(x,y) E K(x,x)].$
\par Let $p =1$. Therefore, in accordance with Formulas
$(10,11)$ and 3.2$(1,2)$ the equality
\par $(12)$ $(\mbox{}_3\sigma _t + ~ \mbox{}_1\sigma ^3_x+ ~ 3  ~
\mbox{}_2\sigma _y ~ \mbox{}_1\sigma ^2_x + ~ 3 ~ \mbox{}_2\sigma
^2_y ~ \mbox{}_1\sigma _x + ~ \mbox{}_2\sigma ^3_y) K(x,y)
$\par $ + 6 (\mbox{}_1^1\sigma _x + ~ \mbox{}_2^1\sigma _y) [
K(x,y)(~\mbox{}_1\sigma _x E K(x,x))] - K(x,y) \{ [
\mbox{}_1\sigma _z, ~ \mbox{}_1\sigma _x] EK(x,z)]|_{z=x}
\}$\par $ - [ \mbox{}^1_1\sigma _x, ~ \mbox{}^2_1\sigma _x] [
K(x,y) EK(x,x)] =0$ \\ follows, when the operator $(I-{\sf A}_xE)$ is
invertible.

In view of Theorem 3.3 this implies:

\par {\bf 2.1. Theorem.} {\it If suppositions of Theorem 3.3 and Example 2
are satisfied. Then a solution of PDE 2$(12)$ with $\mbox{}_1\psi _0=\mbox{}_2\psi _0=0$ over the
Cayley-Dickson algebra ${\cal A}_r$ with $2\le r\le 3$ is given by
Formulas $(2-4)$ and 2$(1)$, when $p =1$.}

\par {\bf 2.2. Example. Korteweg-de-Vries' type PDE.} Continuing Example 2
mention that on the diagonal $x=y$ the operators are: $L_{1,x,x}=0$, $L_{2,x,x} = \partial /\partial t + 8 ~ \mbox{}_1\sigma ^3_x$.
Therefore, $[L_{1,x,x},E]=0$ is valid. Let $E$ be independent of the parameter $t$, then
$[\partial /\partial t , E] =0$, since $t\in \bf R$ and $\bf R$ is the center
of the Cayley-Dickson algebra ${\cal A}_r$. To the term $\mbox{}_1\sigma ^3_x$ the cubic
form $(Im ~ w)^3=-|w|^2w$ corresponds, since $\mbox{}_1\psi _0=0$, where $Im ~w = (w-w^*)/2$,
$w=i_1x_1\psi _1+...+i_{2^{r-1}}x_{2^{r-1}}\psi _{2^{r-1}}$, $~ x_j\in \bf R$ for each $j$.
That is for $E=ST_g$ the restriction is $[|w|^2w,E]=0$, where $n=1$ and $B=1$. Geometrically in the real shadow of
$Im ({\cal A}_r)$ such $E=E(x)$ permits any rotations along the axis $J_w$ parallel to $w$ such that $J_w$ crosses the origin
of the coordinate system. Evidently $[w^3,E]=0$ is satisfied if $[w,E]=0$, that is $[\mbox{}_1\sigma _x, E(x)]=0$.
In the latter case and when $n=1$, $\mbox{}_1\sigma
=\mbox{}_2\sigma $, $\mbox{}_1\psi _0=0$ and $\mbox{}_3\sigma
_t=\partial /\partial t_0$ the differentiation of 2$(12)$ with the
operator $\mbox{}_1\sigma _x$ and the restriction on the diagonal $x=y$ provides the PDE
\par $(1)$ $v_t(t,x)+6~\mbox{}_1\sigma _x[v(t,x)Ev(t,x)] +
~\mbox{}_1\sigma _x^3v(t,x)=0$\\
of Korteweg-de-Vries' type, where $v(t,x)=2~\mbox{}_1\sigma _x
K(x,x)$. Particularly there are solutions of PDE $(1)$ which have the symmetry property
$Ev(t,x) = v(t,x)$.

\par {\bf 3. Example. Non-isothermal flow of a non-compressible
Newtonian liquid with a dissipative heating.} Take the pair of PDOs
\par $(1)$ $L_1=\sigma _x+ \sigma _y$ and
\par $(2)$ $L_2= \mbox{}_1\sigma _t +\sigma _x^2 + q \sigma _y\sigma _x
+\sigma _y^2$, \\ where $q\in {\bf R}$ is a real constant, and consider the
integral equation 3.2$(1)$ with $N$ of the form 3.1$(1)$, so that
\par $(3)$ $L_1F(x,y)=0$ and
\par $(4)$ $L_{2,j}F(x,y)=0$ for each $j$, (see also $(4.81)$ and $(4.82)$ in \cite{lucveleq2013}).
Therefore, in $(4.83)$ \cite{lucveleq2013} the term $K$ changes into $EK$ and
due to Proposition 2.5 and Corollary 2.6 we deduce the formula:
\par \par $(5)$ $L_2K(x,y) = I_1+I_2$, \\ where
\par $(6)$ $I_1= p (~\mbox{}_1^2\sigma _t+ \sigma _x^2 + q~\mbox{}^1\sigma
_y~\sigma _x)~\mbox{}_{\sigma }\int_x^{\infty } F(z,y) E
K(x,z)dz$\par $ =  p ( \mbox{}^2\sigma _x^2 +
q~\mbox{}^1\sigma _y~\mbox{}^2\sigma _x) ~\mbox{}_{\sigma
}\int_x^{\infty } F(z,y) EK(x,z)dz $\par $ -  p
\sigma _x[ F(x,y)E K(x,x)] - p ~\mbox{}^2 \sigma _x
[F(x,y)E K(x,z)]|_{z=x} - q p \sigma _y [F(x,y)E K(x,x)]$ and
\par $(7)$ $ ~ -I_2 = p (\mbox{}^1\sigma _z^2 +
q~\mbox{}^1\sigma _y~\mbox{}^1\sigma _z) ~\mbox{}_{\sigma
}\int_x^{\infty } F(z,y) E K(x,z)dz$
\par $ =  p (\mbox{}^2\sigma _z^2 - q ~\mbox{}^1\sigma _y~\mbox{}^2\sigma _z)
~\mbox{}_{\sigma }\int_x^{\infty } F(z,y) E K(x,z)dz $\par
$ -  p ~\mbox{}^1\sigma _x[ F(x,y) E K(x,x)] + p
~\mbox{}^2 \sigma _z [ F(x,y) E K(x,z)]|_{z=x} - q p
\sigma _y [ F(x,y) E K(x,x)]$. \par Then in $(4.88,4.89)$ \cite{lucveleq2013}
$K$ changes into $EK$ as well and $(4.90)$ takes the form:
\par $(8)$ $L_2 K(x,y) =  p (\mbox{}^2_1\sigma _t+ ~\mbox{}^2\sigma _x^2
+ q ~~\mbox{}^2\sigma _z ~\mbox{}^2\sigma _x  + (q-1)
~~\mbox{}^2\sigma _z^2 ) ~\mbox{}_{\sigma }\int_x^{\infty }
F(z,y) E K(x,z)dz $
\par $- (q+2) p ~~\mbox{}^2\sigma _x [F(x,y) E K(x,x)] +q p
~~\mbox{}^2\sigma _x[F(x,y) E K(x,z)]|_{z=x} +q p~
~\mbox{}^2\sigma _z[ F(x,y) E K(x,z)]|_{z=x}$.
\par Take $q=2$. If the cases of $F$ and $K$ are the same as in
Example 4.6 \cite{lucveleq2013}, $[L_1,E]=0$ and $[L_2,E]=0$ and
when conditions of Theorem 3.3 are fulfilled, the equality follows:
\par $(9)$ $(\mbox{}_1\sigma _t + \sigma _x^2
+ 2 ~\mbox{}^2\sigma _y ~\mbox{}^2\sigma _x  + \sigma _y^2 )
K(x,y) = - 2  p K(x,y) [\sigma _x E K(x,x)] $, \\ where $K$ depends on the parameter $t$.
\par Let $g(x,t)=K(x,x)$, then on the
diagonal $x=y$ this implies the PDE:
\par $(10)$ $(\mbox{}_1\sigma _t + \sigma _x^2) g(x,t) =
- 2 p  g(x,t)[\sigma _x Eg(x,t)]$.
\par Then Equality $(3)$ and Proposition 2.5 imply that $K(x,y) = K(\frac{(\psi , x-y>}{2}$
and the condition $[L_1,E]K=0$ is fulfilled, when $E(x,y) =  E(\frac{(\psi , x-y>}{2}$.
\par Mention that  $L_{2,x,x} =  \mbox{}_1\sigma _t + 4 \sigma _x^2 $ on the diagonal $x=y$.
Taking $E$ independent of $t$, the condition $[L_{2,x,x},E(x,x)]=0$ means that
$[L_{2,x,x},E(0)]=0$. Thus in the real shadow of $Im ({\cal A}_r)$ this $E(0)$ induces any
element of the orthogonal group $O(2^r-1)$. Then more general PDE $(10)$ can be applied to non-isothermal flow
of a non-compressible Newtonian liquid with a dissipative heating as
in \cite{lucveleq2013}. This also provides symmetry properties of $g(t,x)$, particularly,
when a solution satisfies the condition $Eg=g$.

\par {\bf 3.1 Theorem.} {\it Suppose that conditions of Theorem 3.3 and
Example 3 are satisfied, then PDE $(9)$ over the Cayley-Dickson algebra
${\cal A}_r$ with $2\le r\le 3$ has a solution given by Formulas $(3,4)$,
3.1$(1)$ and 3.2$(1)$, where PDOs $L_1$ and $L_2$ are given by
2.1$(1,2)$, $F\in Mat_{n\times n}({\bf R})$ and $K\in
Mat_{n\times n}({\cal A}_r)$, $n\in {\bf N}$ for $r=2$, $n=1$ for $r=3$.}

\par {\bf Conclusion.} In the paper new integrable PDEs were found with the help
of non-commutative integration over octonions and Cayley-Dickson
algebras. It enlarges possibilities of previous approaches based on
real and complex numbers, because each PDE over them can be
reformulated over octonions and new types  of PDEs can be
encompassed. There is the vast general research theme on
integrability of differential equations and PDEs over real and
complex numbers basing on Lie groups and algebras. It is interesting
to develop this theme further and investigate integrable PDEs using
nonassociative analogs of Lie groups and algebras over octonions and
Cayley-Dickson algebras. It is planned to be continued in a next
paper.

\end{document}